\documentclass{article}
\usepackage[normalem]{ulem}
\usepackage{amsthm}
\usepackage{mathtools}
\usepackage{thmtools}

\declaretheoremstyle[headfont=\normalfont]{normalhead}
\usepackage{color}
\usepackage{graphicx}
\usepackage{morefloats}
\usepackage{hyperref}
\usepackage{mathrsfs}
\hypersetup{
    colorlinks=false,
    linkcolor=cyan,
    filecolor=cyan,      
    urlcolor=cyan,
}
\usepackage{ amssymb }
\usepackage{textcomp}
\usepackage{mathabx }
\usepackage{url}
\usepackage{float}
\usepackage{biblatex}
\newtheoremstyle{mydef}
{\topsep}{\topsep}%
{}{}%
{\itshape}{}
{\newline}
{%
  \rule{\textwidth}{0.0pt}\\*%
  \thmname{#1}~\thmnumber{#2}\thmnote{\-\ #3}.\\*[-1.5ex]%
  \rule{\textwidth}{0.0pt}}%

\begin{document}
\newtheorem{conjecture}{Conjecture}
\newtheorem{theorem}{Theorem}
\newtheorem*{theorem-non}{Theorem}
\newtheorem{remark}{Remark}
\newtheorem{proposal}{Proposal}
\newtheorem{proposition}{Proposition}
\newtheorem*{proposition-non}{Proposition}
\newtheorem{lemma}{Lemma}
\newtheorem{corollary}{Corollary}
\newtheorem{observation}{Observation}
\newtheorem{definition}{Definition}
\newtheorem*{question-non}{Question}
\newtheorem*{corollary-non}{Corollary}
\author{Barry Brent}

\date{15 September 21}

\title{Polynomial
interpolation of modular forms for Hecke groups
}
\maketitle
\begin{abstract}
\hskip -.2in
For $m = 3, 4, ...$,
let $\lambda_m = 2 \cos \pi/m$ and
let $J_m (m = 3, 4, ...$) be triangle
functions for the Hecke groups $G(\lambda_m)$
with Fourier
expansions
$
J_m(\tau) = \sum_{n=-1}^\infty a_n(m)q_m^n,
$
where
$q_m(\tau) = \exp 2 \pi i \tau/\lambda_m$.
(When normalized appropriately,
$J_3$ becomes Klein's $j$-invariant
$
j(\tau) = 1/e^{2 \pi i \tau} + 744 + ....
$)
For  $n = -1, 0, 1, 2$ and $3$,
Raleigh gave polynomials $P_n(x)$ such that 
$a_{-1}(m)^n q_m^{2n+2} a_n(m) = P_n(m)$
for $m = 3, 4, ...,$
and conjectured that similar 
relations hold for all positive integers
$n$.
This was proved by Akiyama.
We apply work of Hecke
to study experimentally
similar polynomial interpolations
of the $J_m$ Fourier coefficents
and the Fourier coefficients of
other, positive weight, 
modular forms for $G(\lambda_m)$.
We connect these polynomials 
(again, only empirically)
with variants of Dedekind's eta function, 
with the Fourier expansions of some standard
Hauptmoduln, and, in the case of
analogues of Eisenstein series
for $SL(2,\mathbb{Z})$, with certain
divisor sums.
\end{abstract}
\section{Introduction}
\subsection{An example}
Here is an example of a sequence
$\{P_n(x)\}$
from $\mathbb{Q}[x]$ and a corresponding sequence
of modular forms $\{f_m\}$  having the
relationship
we examine in this article.
Let $T_m$ :=
the cyclic subgroup of $SL(2,\mathbb{R})$
generated by
$$
\left(
\begin{array}{cc}
1 & 2\pi/m\\
0 & 1\\
\end{array}
\right),
$$
let $f_m(x): = \sin(m x)$,
and let
$$
Q_n(x): = (-1)^{(n-1)/2} x^n/n!.
$$
Furthermore let 
$P_n(x) = Q_n(x)$
if $n$ is odd and $P_n(x) = 0$
if $n$ is even.
Members of $SL(2,\mathbb{R})$
act on $\mathbb{R}$ as follows.
If 
$$
M =
\left(
\begin{array}{cc}
a & b\\
c & d\\
\end{array}
\right)
$$
and $x$ is real, we set
$$M(x): = \frac{ax+b}{cx+d}.$$
Thus the $T_m$ act on 
$\mathbb{R}$ by translation.
From the periodicity and Taylor series
of sine,
we know that $f_m(x)$ is invariant
(weight-$0$ modular)
with respect to
to the action
of the  $T_m$ and equal to
$\sum_{n=0}^{\infty} P_n(m)x^n$.
We say that the elements of 
$\{f_m\}$
are interpolated by the 
sequence of polynomials
 $\{P_n(x)\}_{n = 0, 1, ...}$.
\subsection{First sketch of the background.}
Let $\mathbb{Z}, \mathbb{Q}, 
\mathbb{C}$ and $\mathbb{H}$ 
denote, respectively, the set of rational integers,
the set of rational
numbers, the set  of 
complex numbers, and the set of
complex numbers
with positive imaginary parts. 
(We will reserve 
the letter $\tau$ for elements of the
upper half-plane, and $z$
for generic complex numbers.)
We write
$\mathbb{H}^* = \mathbb{H} \cup \mathbb{Q} \cup 
\{i \infty\}$, and we equip $\mathbb{H}^*$
with the Poincar{\'e} metric. Figures
$T$ made by three geodesics
of $\mathbb{H}^*$
are called hyperbolic or
circular-arc triangles.
Let $\lambda_m = 2 \cos \pi/m$.
For $m = 3, 4, ...$,
we define the Hecke group
$G(\lambda_m)$
as the discrete group generated by 
the maps $z \rightarrow -1/z$
and $ z \rightarrow z + \lambda_m$.
The full modular group $SL(2,\mathbb{Z})$
is identical to $G(\lambda_3)$.
\newline  \newline \noindent
To define modular forms for
the Hecke groups, we preview
a definition from Berndt  
\cite{berndt2008hecke}, 
which we will quote again
in a later section.
(We depart occasionally from Berndt's
choices of variable to avoid 
clashes with some of our other
notation.)
\newline \newline \noindent
We say that $f$ belongs to the space
$M(\lambda, k, \gamma)$ if
\begin{enumerate}
    \item 
    $$ f(\tau) = \sum_{n = 0}^{\infty}
    a_n e^{2 \pi i n \tau/\lambda},
    $$
    where $\lambda > 0$ and 
    $\tau \in \mathbb{H}$, and
    \item  $f(-1/\tau) = \gamma(\tau/i)^k f(\tau)$,
    where $k > 0$ and $\gamma = \pm 1$.
\end{enumerate}
We say that $f$ belongs to the space
$M_0(\lambda,k, \gamma)$ if $f$ satisfies conditions
1 and 2 and if $a_n = O(n^c)$
for some real number $c$, as $n$ tends to $\infty$.
\newline \newline \noindent
Members of $M(\lambda, k, \gamma)$
are known as modular forms for 
$G(\lambda)$ of weight $k$.
Condition 1 tells us that  they
are invariant under translations
$\tau \mapsto \tau + \lambda$.
Next we preview Berndt's definition
of cusp forms for Hecke groups.
If $f\in M(\lambda, k, \gamma)$
and $f(i \infty) = 0$,
then we call $f$ a cusp form of
weight $k$ and multiplier  $\gamma$
with respect to $G(\lambda)$. 
For cusp forms, the constant terms 
of condition 1 vanish.
We denote by $C(\lambda,k,  \gamma)$ the
vector space of  all cusp forms
of this kind. 
\newline \newline \noindent
For our purposes,
Schwarz triangles $T$ are 
hyperbolic triangles in $\mathbb{H}^*$
with certain restrictions
on the angles at the vertices.
From a Euclidean point of view, 
their sides are
vertical rays, 
segments of vertical rays, 
semicircles 
orthogonal to the real axis
and meeting it
at points $(r,0)$
with $r$ rational,
or arcs of such semicircles.
We choose
$\lambda, \mu$ and $\nu$,
all non-negative, such that
$\lambda + \mu + \nu <1$;
then the angles of $T$
are $\lambda \pi, \mu \pi$,
and $\nu \pi$.
By reflecting $T$ 
across one of its edges, we
get another Schwarz triangle. 
The reflection between
two triangles in $\mathbb{H^*}$
is effected by 
a M\"obius transformation,
so the orbit of $T$
under repeated reflections
is associated to a 
 collection of 
 M\"obius transformations.
The group 
 generated by 
these transformations is a
triangle group.
By the Riemann Mapping Theorem
there is a conformal, onto 
map
$\phi: T \mapsto \mathbb{H}^*$
called a triangle function.
\newline \newline \noindent
Hecke groups 
are triangle groups 
$H$ that act 
properly discontinuously 
on $\mathbb{H}$.
\footnote{\cite{hecke1936bestimmung}}
This means that for compact 
$K \subset \mathbb{H}$, the set
$\{\mu \in H$ s.t. \hskip -.05in
$K \cap \mu(K) \neq \emptyset\}$
is finite.
Recall that 
$G(\lambda_m)$ is the Hecke group
generated by the maps 
$z \mapsto -1/z$
and $z \mapsto z + \lambda_m$.
Hecke established in \cite{hecke1936bestimmung}
that $G(\lambda_m)$ has the structure
of a free product 
of cyclic groups
$C_2 * C_m$,
generalizing the relation
\cite{serre1970course, cangul1996group}
$SL(2,\mathbb{Z}) = C_2 * C_3$.
\newline \newline \noindent
Let $\rho = -\exp(-\pi i/m) = 
-\cos(\pi/m) + i\sin(\pi/m)$, and
let $T_m \subset \mathbb{H}^*$ 
denote the 
hyperbolic triangle with vertices
$\rho, i$, and $i\infty$. 
The corresponding angles are 
$\pi/m, \pi/2$ and $0$
respectively. Let 
$\phi_{\lambda_m}$ be a 
triangle function 
for $T_m$. The
function $\phi_{\lambda_m}$ 
has a pole at $i\infty$
and period $\lambda_m$.
For $P, Q \in \mathbb{H}^*$,
let us us write $P \equiv_H Q$ when 
$\mu \in H$ and $Q = \mu(P)$.
Then $\phi_{\lambda_m}$ extends to 
a function
$J_m: \mathbb{H}^* \rightarrow \mathbb{H}^*$ by
declaring that
$J_m(P) = J_m(Q)$
if and only if 
$P \equiv_H Q$.
$J_m$ is a modular function
for $G(\lambda_m)$.
\newline \newline \noindent
Schwarz, Lehner, Raleigh and others
studied Schwarz triangle functions, 
which map  hyperbolic triangles $T$ in the 
 extended upper half $z$-plane  onto
 the extended upper half 
 $w$-plane.\footnote{\cite{schwarz1873ueber},
\cite{lehner1954note},\cite{raleigh1962fourier}}
 For certain $T = T_m$, 
a triangle function 
$\phi_{\lambda_m}: T \to \mathbb{H}^*$
 extends to a map
$J_m: \mathbb{H}^* \rightarrow \mathbb{H}^*$
 invariant under
modular transformations
from $G(\lambda_m)$.
Suitably normalized, the $J_m$
become analogues 
$j_m$ of the normalized
Klein's modular invariant 
$$
j(\tau) = 1/q + 744 + 196884 q + ...
$$
where $q = q(\tau) = \exp (2 \pi i \tau)$ 
 and 
 $j_3(\tau) = j(\tau)$.\footnote{For $j$, see 
 \cite{serre1970course}, Chapter VII, equation (23).}
 The $j_m$ are studied in conjecture 1 below.
 \newline \newline \noindent
With $\lambda_m = 2 \cos \pi/m$ and 
$q_m(\tau) = \exp (2 \pi i \tau/\lambda_m)$,
the original $J_m$ have Fourier series
$J_m(\tau)= \sum_{n \geq -1} a_n(m) q_m(\tau)^n$.
For  $n = -1, 0, 1, 2$ and $3$,
Raleigh 
gave polynomials $P_n(x)$ such that 
$a_{-1}(m)^n q_m^{2n+2} a_n(m) = P_n(m)$
for $m = 3, 4, ...,$
and conjectured that similar 
relations hold for all positive integers
$n$. \footnote{\cite{raleigh1962fourier}}
Akiyama 
proved this conjecture
 in the passage 
after his (Akiyama's) equation (6)
\cite{akiyama1992note}.
\newline  \newline \noindent Hecke 
built families of modular forms $f_m$
 for $G(\lambda_m)$
sharing particular properties.
\footnote{\cite{berndt2008hecke}, 
\cite{hecke1936bestimmung}}
Earlier authors, whose work we will 
also describe,
had already built modular functions
(meromorphic functions
invariant under the
action of $G(\lambda)$,
thus, of weight zero)
from triangle functions.
\subsection{Plan of the article.}
The plan of the article is as follows.
(1) An elaboration of the preceding discussion
to establish a basis for the  code 
in our experiments.\footnote{We 
have documention in
the data repository \cite{test3}. 
\it Mathematica \rm
notebook names
end in the suffix ``.nb'', and
\it SageMath \rm notebook names
end in the suffix ``.ipynb''.
Numerical data files named in the notebooks
is stored in the folder ``data''
on \cite{test3}. 
A green ``Code'' button on the top page of 
the repository contains a drop-down menu
with a download option.
A \it Mathematica \rm notebook 
(``mf25.nb'') in the repository
is a searchable library 
of functions that may not be defined
explicitly within our other notebooks.
We used 
\it SageMath \rm release 9.1.}
(2) Conjectures on
polynomials in $\mathbb{Q}[x]$
interpolating the coefficients in Fourier expansions
 of triangle functions for $G(\lambda_m)$.
(3) A survey of Hecke's theory
of modular forms for $G(\lambda_m)$,
especially, the construction
of modular forms 
from modular functions.
(4) Several conjectures
about polynomials 
in $\mathbb{Q}[x]$
interpolating the 
 coefficients in Fourier expansions
 of Hecke modular
forms on $G(\lambda_m)$.
(5) Several data plots and tables.
Tables at the end of the article
focus on the triangle functions, 
since they are the basis of
our construction of positive-weight
modular forms, but
more extensive collections of 
plots and tables are available
within the \it Sagemath \rm and
\it Mathematica \rm notebooks 
on \cite{test3}.
\subsection{Methods.}
Our conjectures are based on
numerical experiments;
here is a little  more detail
on the way we arrive at them.
We begin with a list of modular functions or
modular forms $f_m$
for $G(\lambda_m), m = 3, 4, ...$ sharing
certain properties picked out by
Hecke's theory.
Then we make tables 
of polynomials 
$Q_n(x)$ generated by
Lagrangian interpolation
from the values of the coefficient $k_m(n)$ 
in Fourier expansions
$f_m = \sum_n k_m(n) X_m^n$,
where $X_m$ is a variable related
to $q_m(\tau)$. Thus we are
seeking $Q_n(x)$ such that
\begin{equation}
Q_n(m) = k_m(n)
\end{equation}
for $m = 3, 4, ....$.  
If the degrees of the $Q_n(x)$ we obtain
are linear in $n$, we take
this to be evidence that the $Q_n(x)$
do satisfy equation (1)
 for all integers $m$ greater than two. 
 (Typically,
the alternative outcome
is that the degree of every polynomial
$Q_n(x)$ that we 
generate in a given table is equal
to the size of the data set we are 
trying to interpolate.)
\subsection{Work of Lehner, Raleigh and Leo.}
The earliest 
computer code
we located 
for calculating 
Fourier 
expansions of 
triangle functions for Hecke groups
is that of Leo\footnote{\cite{leo2008fourier}}; 
it is based on Lehner's 
construction. Leo also
calculates the Fourier coefficients
of weight 4 and weight 6
Hecke-analogues of classical Eisenstein
series in Chapter 4 of \cite{leo2008fourier}.
Our code for triangle functions,
which is based on Leo's,
comes from the papers of 
Lehner and Raleigh.
J. Jermann's package
is also concerned with modular forms
of triangle groups for Hecke groups,
but we did not make use of it.
\footnote{\cite{sagemath}}
\subsection{Disclaimer.}
The article describes experiments
and states conjectures. It contains
no theorems except ones
that we quote from
the existing literature.
\section{A glossary}
Some special functions in this list 
are related; different notations for 
similar objects
are used by Lehner and Raleigh,
and we included all of them.
\begin{enumerate}
\item The digamma function
$\psi(z):=\Gamma'(z)/\Gamma(z)$.
\item \footnote{\cite{caratheodory2}, p. 130, 
equation 370.8}
The Schwarzian derivative
\begin{equation}
\{w,z\} = 
\frac{2 w' w''' - 3w''^2}{2w'^2}
\end{equation}
for $w = w(z)$.
(In section 3 below,
we discuss Caratheodory's
presentation of a
well-known theorem of Schwarz;
Carath{\'e}odory writes the left side
of our equation (4)  as 
``$\{w,z\} = \frac{w^{'}w^{'''} - 
3 w^{''2}}{w^{'2}} = ...$''
\footnote{\cite{caratheodory2}, \textsection 374}, 
but we infer that the Schwarzian derivative $\{w,z\}$ 
is intended
from the automorphy property of 
that theorem's clause 2.)
\item The Pochhammer  symbol
$$
(a)^0: = 1 \text{ and, for }
n \geq 1,
(a)^n:=a(a+1)...(a+n-1) = \Gamma(a+n)/\Gamma(a).
$$
\item 
\footnote{\cite{caratheodory2}, 
p. 138, equation 377.3}
The function $c_{\nu}$ given by 
$$
c_{\nu} = c_{\nu}(\alpha, \beta, \gamma) :=  
\frac{(\alpha)^{\nu} (\beta)^{\nu}}
{{\nu}!(\gamma)^{\nu}}, {\nu} \geq 0.
$$
To facilitate comparison with
Raleigh's  
equation ($9^1$) 
\footnote{\cite{raleigh1962fourier}}
, we remark that 
\begin{equation}
c_{\nu} = 
\frac{\Gamma(\alpha + {\nu})}{\Gamma(\alpha)}
\cdot \frac{\Gamma(\beta + {\nu})}{\Gamma(\beta)}
\cdot \frac{\Gamma(1)}{\Gamma(1 + {\nu})}
\cdot \frac{\Gamma(\gamma)}{\Gamma(\gamma + {\nu})}.
\end{equation}
In the terms of this article's Theorem 1 
below, Raleigh is treating the
case $\lambda = 0$, for which 
(equation (7) below) $\gamma = 1$
and the expression on the right side of (3)
becomes, as in Raleigh,
$$\frac{\Gamma(\alpha + {\nu})\Gamma(\beta + {\nu})}
{\Gamma(\alpha)\Gamma(\beta)({\nu}!)^2}.$$
\item 
The function $e_{\nu}$ given by
\footnote{\cite{raleigh1962fourier}, equation $9^1$}
$$
e_{\nu} = e_{\nu}(\alpha, \beta) :=  
\sum_{p = 0}^{{\nu} - 1}
\left ( \frac 1{\alpha + p} + 
\frac 1{\beta + p} 
- \frac 2{1 + p} \right ).
$$
Here, we are dealing with the 
same ambiguity present in the
definition of $c_{\nu}$:
this is a specialization to the
case $\gamma = 1$ of the
$e_{\nu}$ 
for ${\nu} \geq 1$
given by 
\footnote{\cite{caratheodory2}, p. 153,
equation 387.5}
$$
e_{\nu} = e_{\nu}(\alpha, \beta, \gamma) :=  
\sum_{p = 0}^{n - 1}
\left ( \frac 1{\alpha + p} + 
\frac 1{\beta + p} 
- \frac 2{\gamma + p} \right ).
$$ Unless it
is explicitly indicated to be
otherwise, we intend the former (Raleigh's)
definition.
\item 
\begin{enumerate}
    \item 
\footnote{\cite{caratheodory2}, p. 138,
equation (377.4)}
    Gauss's hypergeometric series
$$
F(\alpha,\beta,\gamma;z):=
\sum_{{\nu}=0}^{\infty} 
c_{\nu}(\alpha, \beta, \gamma) z^{\nu}.$$
$F$ is occasionally written in 
\cite{caratheodory2} as $\phi_1$ (for
example, on, p. 152.)
    \item 
 \footnote{As
    defined in the first line of
\cite{caratheodory2}, p. 142.}   
$$
F_1(\alpha,\beta,\gamma; z):=
F(\alpha,\beta,\gamma+1;z).
$$
\item 
\footnote{\cite{leo2008fourier}, equation (3.5)}
Alternatively,
    dropping $\gamma$:
$$
F_1(\alpha,\beta;z):= 
\sum_{\nu=1}^{\infty}
\frac{(\alpha)_k(\beta)_{\nu}}{(\nu!)^2}
e_{\nu}(\alpha, \beta).$$
It is in the latter form, defined
more cryptically in \cite{lehner1954note}, p. 244,
that we will use $F_1$; to
establish his series for the triangle functions,
which we will apply below,
Lehner uses this definition of
$F_1$, as well as certain 
theorems from  Fricke.
\footnote{\cite{fricke1922elliptischen}}
Referring to item 4, we see that
$$
F_1(\alpha,\beta;z)= 
\sum_{\nu=1}^{\infty}
c_{\nu}(\alpha,\beta,1)
e_{\nu}(\alpha, \beta).$$
We will derive another form of $F_1(\alpha,\beta;z)$
in item 7.
\end{enumerate}
\item With 
$F = F(\alpha, \beta, \gamma;z)$,  
a special function
$$F^*(\alpha, \beta, \gamma;z): = 
\frac{\partial F}{\partial \alpha} + 
\frac{\partial F}{\partial \beta}+ 
2\frac{\partial F}{\partial \gamma};$$
$F^*$ may be written
\footnote{\cite{caratheodory2}, equation (387.4) on p. 153}
$$
F^*(\alpha, \beta, \gamma;z) =
\sum_{{\nu} = 1}^{\infty}
c_{\nu} (\alpha, \beta, \gamma) 
e_{\nu} (\alpha, \beta, \gamma) z^{\nu}.
$$
It follows that $F^*(\alpha, \beta, 1;z) = F_1(\alpha,\beta;z)$.
\item 
\footnote{\cite{caratheodory2},
p. 152, equations 386.2 and 386.3}
A special  function
$\phi_2^*(z)$ is defined as a certain limit
but is immediately 
reduced to 
$$\phi_2^*(z) =  
F(\alpha, \beta, 1;z) \log z
+ F^*(\alpha, \beta, 1;z).$$
\item 
\footnote{\cite{OEIS11}}
The set $\mathscr{Q} 
=\{2,5,6,8,10,11,14,15,17,18,20,22,23,...\}$
of positive integers not
represented by the quadratic form
$x^2+xy+y^2$ . B. Cloitre
asserts on the cited page
that  $\mathscr{Q}$ is also the set of
non-negative integers $n$
such that
$\delta(n)$ is non-zero, where
$\eta$ is Dedekind's eta function and
$\sum_n \delta(n)x^n = \eta(x^3)/\eta(x)^3.$
\item The McKay-Thompson series of class 4A,
$\{1,24, 276, 2048, ...\}$, which is
the sequence of coefficients 
in the $q$-series of a certain
hauptmodul discussed in \cite{doi}.
We identified it with the
sequence $\{\phi_n\}$  of our
conjecture 1 after finding it 
in on \cite{OEIS2}.
\item As usual, the 
cardinality of a finite set $S$
is written $\# S$, the
$n^{th}$ prime number is
denoted by $p_n$,
the number of primes less than or equal to $x$
is written $\pi(x)$,
and $\sigma_k(n) :=\sum_{0<d|n} d^k$.
\end{enumerate}
\section{Calculation of Schwarz's inverse 
triangle function}
Schwarz
proved 
\begin{theorem}
\footnote{\cite{caratheodory2}, \textsection 374}
\begin{enumerate}
\item
Let the half-plane $\Im z > 0$ be mapped
conformally onto an arbitrary 
circular-arc triangle
whose angles at its vertices $A, B$, and $C$ 
are $\pi \lambda, \pi \mu$, and $\pi \nu$, and
let the vertices $A, B, C$ be the images of
the  points $z = 0, 1, \infty$, respectively.
Then the mapping function $w(z)$ must be 
a solution of the third-order differential
equation
\begin{equation}
\{w,z\} =
\frac{1-\lambda^2}{2 z^2} +
\frac{1-\mu^2}{2(1-z^2)} +
\frac{1-\lambda^2-\mu^2+\nu^2}
{2z(1-z)}.
\end{equation}
\item
If $w_0(z)$ is any solution
of equation (4) that satisfies
$w'_0(z) \neq 0$ at all interior points
of the half-plane, then the function
$$
w(z) = \frac {aw_0(z) + b}{cw_0(z) + d}
\hskip 1in (ad - bc \neq 0)
$$
is likewise a solution of equation 3.
\item
Also, every solution of equation (4) that is regular
and non-constant in the half-plane $\Im z > 0$
represents a mapping of this half-plane onto a
circular-arc triangle with angles $\pi \lambda,
\pi \mu,$ and $\pi \nu$.
\end{enumerate}
\end{theorem}
\noindent
(In Carath{\'e}odory's lexicon,
a regular function is one that is
differentiable 
on an open connected set.)
\footnote{\cite{caratheodory1} p. 124}
\newline \newline \noindent
Let us write
\begin{equation}
\alpha = \frac 12
(1 - \lambda - \mu + \nu),
\end{equation}
\begin{equation}
\beta = \frac 12(1-\lambda - \mu - \nu),
\end{equation}
and
\begin{equation}
\gamma = 1 - \lambda.
\end{equation}
The solutions $w$ of equation (4)
are inverse to triangle functions;
they are quotients of 
arbitrary solutions of
\begin{equation}
u'' + p(z)u' +q(z)u = 0
\end{equation}
when 
 $$
 p = \frac {1-\lambda}z-\frac {1-\mu}{1-z}
 $$
 and
  \footnote{\cite{caratheodory2}, 
p. 136, equation (376.4)}
 $$
 q = - \frac {\alpha \beta}{z(1-z)}.
 $$
Equation (8) reduces
to the 
hypergeometric differential
equation 
\footnote{\cite{caratheodory2}, p. 137,
equations 376.5-7}
\begin{equation}
z(1-z)u'' +(\gamma - (\alpha + \beta +1)z)u'
- \alpha \beta u = 0.
\end{equation}
As long as $\gamma$ is not a
non-positive integer, 
$u=F(\alpha,\beta,\gamma;z)$
is a solution of equation (9); it is 
the only solution
regular at $z = 0$, and it
satisfies 
\footnote{Final paragraph of 
\cite{caratheodory2}, \textsection 377, 
p. 138.}
$F(\alpha,\beta,\gamma;0) = 1$.
 \newline \newline \noindent
 In  \cite{caratheodory2}, we find that
 \footnote{\textsection \textsection 386-388 
(pp. 151-155)}
when $\gamma = 1$
and $\lambda = 0$,
 another, linearly independent,
solution of equation (8) is $\phi_2^*(z)$.
Section 394, pp. 165 - 167 
 of \cite{caratheodory2} is 
devoted to the
case $\lambda = 0$. 
There we find that the mapping 
function $w$ of Theorem 1 
satisfies 
\footnote{\cite{caratheodory2}
p. 166,  equation 394.4}
\begin{equation}
 w = \frac 1{\pi i}
 \left [ \frac{\phi_2^*}{\phi_1} -
\left (2 \psi(1) - 
\psi(1 - \alpha)
- \psi(1-\beta) \right ) \right ] + 
i \frac{ \sin \pi \mu}
{\cos \pi \mu + \cos \pi \nu}.
\end{equation}
\section[]{Inversion of Schwarz's 
inverse triangle function}
Following Lehner and Raleigh,
we consider the Schwarz triangle
$T_m$ with vertices at 
$\rho = -\exp(-\pi i/m), i$,
and $i\infty$. 
In terms of
Theorem 1,
$T_m$ has
$\lambda = 0$
(an angle $0$ at 
the vertex $i\infty$), $\mu = 1/2$
(an angle $\pi/2$ at $i$),
and $\nu = 1/m$ (an angle $\pi/m$
at $\rho$.)
 In this situation, $\gamma = 1$. 
 \newline \newline \noindent
 Let $J_m$ be automorphic for 
 $G(\lambda_m)$ with
 $J_m(\rho) = 0, J_m(i) = 1$,
 and $J_m(i \infty) = \infty$.
 In terms of Theorem 1, $w$
 and $J_m$ are inverse functions.
 We are going to write
 down the Fourier expansion
 $\sum_{n = -1}^{\infty} a_n q_m(\tau)^n$
 of $J_m$.
 \newline \newline \noindent
By clause 2 of Theorem 1, if $w$ satisfies 
equations (4) and (10), so does
$\tau = \thinspace \tau(z) = \lambda_m w(z)/2$, 
and therefore
$$
2\pi i\tau /\lambda_m =
 \frac{\phi_2^*}{\phi_1} - \left (2 \psi(1) - 
\psi(1 - \alpha)
- \psi(1-\beta) \right ) - 
\pi  \sec(\pi/m).
$$
Let us write
$\log A_m = -2 \psi(1) + \psi(1 - \alpha)
+ \psi(1 - \beta)  - \pi \sec(\pi/m).$ 
In general, $A_m = a_{-1}(m)$.
\footnote{\cite{raleigh1962fourier},
for example two lines below equation (13).}
Recalling the definitions of 
$\phi_1$ and $\phi_2^*$
from our glossary items 6 and 8, 
we find (abbreviating $J_m(\tau)$ as $J_m$) that 
\begin{equation}
2\pi i\tau /\lambda_m 
  = - \log J_m
+ \frac {F^*(\alpha, \beta, 1;1/J_m)}
{F(\alpha, \beta, 1;1/J_m)} +  
\log A_m.
\end{equation}
Equation (11) is equation (6) of 
\cite{raleigh1962fourier},
but Raleigh  suppresses the 
subscripts.
He also writes 
$\exp 2 \pi i \tau/\lambda_m$ as $x_m$, 
so that (in our earlier notation)
$x_m = q_{_m}(\tau)$.
\newline \newline \noindent
In Raleigh's notation, after taking
exponentials,
\begin{equation}
 x_m/A_m = 
\frac 1{J_m} \exp 
\frac {F^*(\alpha, \beta, 1;1/J_m)}
{F(\alpha, \beta, 1;1/J_m)},
\end{equation}
the right side of
which has a power series in $J_m$ with
rational coefficients. 
Writing $X_m = x_m/A_m$
we can regard $X_m = X_m(J_m)$ as a power series
in $J_m$ with rational coefficients.
Following \cite{lehner1954note} and
 \cite{raleigh1962fourier},
we inverted this power series to obtain
one for the modular function $J_m$,
also with rational coefficients.
By construction,
the Fourier expansion of $J_m$ in
$X_m$ is normalized so that the coefficient
of $1/X_m$ is $1$. 
\footnote{\cite{raleigh1962fourier}, equation (12).}
Let $\mathscr{I}$
be a formal operation
taking a power series $\sigma(v)$
to its inverse; that is,
if $u=\sigma(v)$ then 
$v = \mathscr{I}(\sigma)(u)$.
Let  $Y_m(J)$ be a power series
such that
$$ Y_m(J_m) = J_m \exp 
\frac {F^*(\alpha, \beta, 1;J_m)}
{F(\alpha, \beta, 1;J_m)} = 
X_m \left (1/J_m \right )$$
and hence
$$
Y_m(1/J_m) =
\frac 1{J_m} \exp _m
\frac {F^*(\alpha, \beta, 1;1/J_m)}
{F(\alpha, \beta, 1;1/J_m)} = X_m(J_m),$$
so that $\mathscr{I}(Y_m)(X_m(J)) = 1/J_m$
and, therefore,
 $J_m = 1/\mathscr{I}(Y_m)(X_m)$.
\begin{remark}
We noticed several typos in
\cite{raleigh1962fourier}.
Four of Raleigh's equations---(I), (10), and
the two equations on p. 109
that begin ``$a_{-1}(q) = ...$''
(where Raleigh's $q$ is our $m$)---are 
pairwise contradictory.
From the second paragraph
on Raleigh's p. 110, we expect that 
$A_3 = a_{-1}(3) = 1/1728, A_4 = a_{-1}(4) = 1/256$, and 
$A_6 = a_{-1}(6) = 1/108$. These values
are consistent with Raleigh's equation 
(10), but not with the others.
We infer that all of them except (10) are incorrect.
Thus, following Raleigh by
writing  $\psi$ for the digamma function,
$\alpha(m)$ for $(1/2-1/m)/2$, and 
$\beta(m)$ for $(1/2+1/m)/2$, 
$$
a_{-1}(m) = \exp\left (-2 \psi(1) + \psi(1-\alpha(m))
+\psi(1 - \beta(m))- \pi \sec(\pi/m)
\right ).
$$
\end{remark}
\section{Raleigh's polynomials for triangle functions}
Let $X_m$ be the variable
from the previous section.
We define some operators on
infinite series in $X_m$.\footnote{The 
substitution
involved appears 
in \cite{leo2008fourier}.}
\begin{definition}
Let $f = \sum_{n=a}^{\infty} k_n X_m^n$, 
where $k_n$ is a rational number for $n = a, a+1, ...,$
and $k_a \neq 0$. 
\begin{enumerate}
\item
Let 
$g = \sum_{n=a}^{\infty} k_n (2^6 m^3 X_m)^n =
\sum_{n=a}^{\infty} \tilde{k}_n X_m^n$ (say). 
Then
$$\overline{f} := g/\tilde{k}_a.$$
\item 
$$
f^*:= \frac 1{k_a} \sum_{n=a}^{\infty} k_n X_m^{n-a}.
$$
\end{enumerate}
\end{definition} \noindent
Recall, from the passage following equation (12)
in the previous section, that
the Fourier expansion of 
$J_m$ in  $X_m$  has the form 
$$
J_m(\tau) = 1/X_m+\sum_{n = 0}^{\infty}a_n(m) X_m^n.
$$
\begin{definition} 
For the present purpose, we regard $J_m$ 
as a Laurent series in $X_m$ and write
$$j_m: = \overline{J_m}.$$
\end{definition} 
\begin{conjecture}
\footnote{Some code for 
 $j_m$ Fourier expansions
appearing in
\it SageMath \rm notebooks
cited below was generated in
`j from scratch.ipynb' \cite{test3}, which
employs a ``dictionary''
(the definitions at the top of the notebook)
distinct from the corresponding dictionaries
in the notebooks where it is
reproduced.}
Let 
the Fourier expansion of 
$j_m(\tau)$ be
$$j_m = 1/X_m +
\sum_{n\geq 0} c_m(n) X_m^n.$$
\begin{enumerate}
\item 
\footnote{Notebook
``conjecture 1.nb'',\cite{test3}.}
For each integer $n$
greater than $-2$,
there exists a polynomial
$C_n(x) \in \mathbb{Q}[x]$
that satisfies the relation
$c_m(n) = C_n(m)$ for $m = 3, 4, ....$
\item 
\footnote{Notebook ``conjecture 1 clause 2.ipynb'',\cite{test3}}
Let $\{\phi_n\}$
be as in item 10 of our glossary. For some
 degree $2n$, irreducible, 
 monic polynomial $\gamma_n(x)$
in $\mathbb{Q}[x]$:
$$
C_n(x) = 
\phi_n \cdot (x-2)(x+2)x^{n+1} \gamma_n(x).
$$
\item 
\footnote{Notebook ``conjecture 1 cause 3.nb'',
\cite{test3}.}
$j_3$ is the
modular function on $SL(2,\mathbb{Z})$
usually denoted $j$.
\item \footnote{Notebooks
``conjecture1clause4.nb'',
``conjecture1clause4d.nb'',
``conjecture 1 clause 4 no2.ipynb'',
``conjecture 1 clause 4 no3.ipynb'',
``conjecture 1 clause 4 no4.ipynb'',
``conjecture 1 clause 4 no5.ipynb'',
and
``conjecture 1 clause 4 no6.ipynb'',
\cite{test3}.
}
The complex roots of 
$\gamma_n(x)$ lie in the disk
with center zero and radius $n/\log (n)$.
\item \footnote{Folder ``conjecture1clause5'',
\cite{test3}.}
Let $G_n$ be the Galois group
of $\gamma_n(x)$ over the rationals. The
size of $G_n$ is $2^n n!$
and (if $n$ is greater than two) 
$G_n$ is isomorphic
to a permutation group on $2n$ elements
$\{e_1, ..., e_{2n}\}$
with three generators: a transposition
$(e_j, e_k)$, a product $(e_j,e_{j'})(e_k, e_{k'})$,
and  a product $\Gamma_1 \Gamma_2$
of disjoint cycles
$\Gamma_1$ and $\Gamma_2$, each of length $n$,
such that $\Gamma_1$ sends $e_j$ to $e_{j'}$
and $\Gamma_2$ sends $e_k$ to $e_{k'}$.
\item \footnote{Notebook ``conjecture 1 clause 6.ipynb'', \cite{test3}.}
 Let $n$ be larger than one and let
$\pi_n$ be the
    set of prime numbers
    dividing the denominator
    of at least one non-zero
    coefficient of $C_n(x)$
    in its unfactored form.
    Then 
    \begin{enumerate}
    \item $\pi_2 = \{3\}$ and $\pi_3$ is empty.
      \item If $\pi_n$ is ordered by size, it contains
        no gaps. That is, if $p$ and $p'$
        are consecutive elements of  $\pi_n$ 
        with $p = p_k$ and $p' = p_j$,
        then $j = k + 1$.
            \item If $n$ is an odd prime other than $3$, 
            then
            $$\pi_n = \{3, ...,k,...,p\}_{k \text{ prime}}$$ where
            $p$ is the greatest prime less than $n$.
            \item If $n$ is composite and
            $n+1$ is prime, then 
            $$\pi_n = \{3...,k, ...,n+1\}_{k \text{ prime}}.$$
            \item If $n$ and $n+1$ are both composite, 
            then  
            $$\pi_n = \{3...,k, ...,p\}_{k \text{ prime}}$$
            where $p$ is the greatest prime less than $n$.
    \end{enumerate}
\end{enumerate}
\end{conjecture} 
 \noindent
Clause 2 implies that, for 
$m$ greater than  or equal to three,
$c_n(m)$ is nonzero.
\noindent
It is already known that,
for all integers $n \geq -1$,
the $n^{th}$ Fourier coefficient of
$j = j_3$, namely
$c(n) = c_n(3)$,
is positive.\footnote{See, for example,
page 199 in \cite{rankinModular}.}
We tested clause 4 in several ways.
We approximated the roots of
the $\gamma_n(x)$  with root-finding routines
and compared their
complex moduli with $n/\log (n)$.
We used the argument principle to count the zeros 
in central disks of radius $n/\log (n)$. 
%
%
We superimposed plots of the roots
of $\gamma_n(x)$ against plots of circles with 
radius $n/\log (n)$
and center at the origin. An example is 
depicted in Figure 1.\footnote{
Notebook ``conjecture1clause4d.nb'',\cite{test3}.}
For clause 5, we computed the Galois 
groups in \it Magma\rm. For clause 6,
some sequences we generated in the
analysis were identified in \cite{OEIS15}
and \cite{OEIS16}.
\begin{conjecture}
\footnote{Relevant documents in \cite{test3}
are notebooks
``conjecture 2.nb'',
``conjecture2no1.ipynb'',
``capital-J make data file1jun21.ipynb''
and associated data files.}
Let the Fourier expansion of $J_m(\tau)$
be
$$
J_m = 
\sum_{n = -1}^{\infty}a_m(n) X_m^n.
$$
\begin{enumerate}
    \item \footnote{For clause 1, see ``conjecture 2.nb'',
    ``conjecture 2 clause 1b.ipynb'',
, and 
``conjecture 2 clause 1b no2.ipynb'',
\cite{test3}.}
We have
    \begin{enumerate}
    \item There exist polynomials $A_n(x)$
    such that
    $A_{-1}(x) \equiv 1$,
    $A_0(x) = 3x^2+4$, 
    $A_1(x) = 69x^4 - 8x^2 - 48$,
    and $A_n(m) = m^{2n+2}a_m(n)$ 
    for $m = 3, 4, ....$.\footnote{The
    first few polynomials in table 10.5
     agree with Raleigh's equation-group
     III
    in \cite{raleigh1962fourier}.}
\item 
\footnote{Notebooks ``conjecture 2.nb'',
``conjecture 2, clause 2.ipynb'',
\cite{test3}.}
Let $C_n(x)$ be as in conjecture 1.
We have:
$$A_n(x) = 2^{-6n-6} x^{-n-1} C_n(x).$$
\end{enumerate}
    \item  \footnote{\cite{test3},
Notebook ``conjecture 1 clause 2 w code 14jun21.ipynb''} 
Let $\pi_n$ be the
    set of prime numbers
    dividing the denominator
    of at least one non-zero
    coefficient of $A_n$.
    Then 
    \begin{enumerate}
    \item $\pi_2 = \{3\}$.
      \item If $\pi_n$ is ordered by size, it contains
        no gaps. That is, if $p$ and $p'$
        are consecutive elements of  $\pi_n$ 
        with $p = p_k$ and $p' = p_j$,
        then $j = k + 1$.
            \item If $n$ is an odd prime, 
            then
            $$\pi_n = \{2, ...,k,...,p\}_{k \text{ prime}}$$ where
            $p$ is the greatest prime less than $n$.
            \item If $n$ is composite and
            $n+1$ is prime, then 
            $$\pi_n = \{2...,k, ...,n+1\}_{k \text{ prime}}.$$
            \item If $n$ and $n+1$ are both composite, 
            then  
            $$\pi_n = \{2...,k, ...,p\}_{k \text{ prime}}$$
            where $p$ is the greatest prime less than $n$.
    \end{enumerate}
\end{enumerate}
\end{conjecture}
\noindent
The existence statement in
clause 1a of conjecture 2 is 
equivalent
up to some changes of variable,
obviously, to Raleigh's 
conjecture.\footnote{(proved in
\cite{akiyama1992note}}
We identified the leading numerical term in clause 1b 
of conjecture 2 after 
looking at \cite{OEIS17}. Clause 2
of conjecture 2
is only a slight refinement of \cite{akiyama1992note},
proposition 2.
\section{Survey of Hecke's theory of modular forms}
When the  
$w$-image of $\mathbb{H}^*$ is
$T_m$, the inverse of $w$ is
$\phi_{\lambda_m}$. 
The extension 
by modularity $J_m$ of $\phi_{\lambda_m}$ to
$\mathbb{H}^*$, is periodic
with period $\lambda_m$ and
maps $\rho$ to $0$, $i$ to $1$, 
and $i\infty$ to 
$\infty$.\footnote{\cite{lehner1954note},
equation (2).}
These mapping properties
allow us, following Berndt's 
exposition
of Hecke,
to construct positive 
weight modular forms for $G(\lambda_m)$
from $J_m$.\footnote{\cite{berndt2008hecke}}
This section describes results of Hecke
that are perhaps most easily accessible
for the classical case $m = 3$ 
in Schoeneberg
 and, for the general case, in 
 Berndt.
 \footnote{\cite{schoeneberg1974},\cite{berndt2008hecke}}
\subsection[]{The case $m=3$.}
By keeping
track
of the weights, zeros and poles of the
constituent factors in the 
numerator and denominator of the
fraction defining
$$ f_{a,b,c} = 
\frac{J^{'a}}
{J^b (J - 1)^c},
$$
Schoeneberg
demonstrates  that $f_{a,b,c}$
is an entire modular form of weight $2a$
for $SL(2,\mathbb{Z})$
 if $a \geq 2, 3c  \leq a, 3b \leq 2a$,  
 $b+c \geq a$ and $a, b, c$ are integers.
(Schoeneberg speaks of ``dimension $-2a$.'')
\footnote{\cite{schoeneberg1974}, 
Theorem 16, p.45}
Thus he is able to 
write down a weight $4$ entire modular form
$E^*_4 = f_{2,1,1}$
for $SL(2,\mathbb{Z})$ with a zero
of order $\frac 13$ at $\rho = e^{2 \pi i/3}$
and a weight $6$ 
entire modular form 
$E^*_6 = f_{3,2,1}$
for $SL(2,\mathbb{Z})$
with a zero of order $\frac 12$
at $i$.
(Schoeneberg writes $G^*_4, G^*_6$.)
It is well known that the 
(vector space) dimension of the
spaces of weight $4$ and $6$ entire
modular forms for $SL(2,\mathbb{Z})$
is equal to one,
so $E^*_4$ and  $E^*_6$ 
may be identified
with the usual weight $4$ 
and weight $6$ 
Eisenstein series,
up to a  normalization. 
Finally, Schoeneberg defines
the weight $12$ cusp form
$\Delta^* = E_4^{*3} - E_6^{*2}$
with a zero of order $1$ at $i \infty$. 
It is a multiple of $\Delta$.
\subsection[]{The case $m \geq 3$.}
We quote statements from Berndt, which is
an exposition of 
Hecke.\footnote{\cite{hecke1938lectures} and
other writings.}
We depart occasionally from Berndt's
choices of variable to avoid 
clashes with our earlier
notation.
\begin{definition}
\footnote{\cite{berndt2008hecke}, Definition 2.2}
We say that $f$ belongs to the space
$M(\lambda, k, \gamma)$ if
\begin{enumerate}
    \item 
    $$ f(\tau) = \sum_{n = 0}^{\infty}
    a_n e^{2 \pi i n \tau/\lambda},
    $$
    where $\lambda > 0$ and 
    $\tau \in \mathbb{H}$, and
    \item  $f(-1/\tau) = \gamma \cdot (\tau/i)^k f(\tau)$,
    where $k > 0$ and $\gamma = \pm 1$.
\end{enumerate}
We say that $f$ belongs to the space
$M_0(\lambda,k, \gamma)$ if $f$ satisfies conditions
1 and 2, and if $a_n = O(n^c)$
for some real number $c$, as $n$ tends to $\infty$.
\end{definition}
\noindent
After defining the notion of a fundamental 
region in the usual
way and defining as $G(\lambda)$ the group 
of linear fractional transformations
generated by
$\tau \mapsto -1/\tau$ and
$\tau \mapsto \tau + \lambda$,
Berndt states
(for $\tau = x + iy$)
\begin{theorem} \footnote{\cite{berndt2008hecke}, Theorem 3.1} 
Let $B(\lambda)= \{\tau \in \mathbb{H}:
x < \lambda/2,  |\tau| > 1\}$. Then if 
$\lambda \geq 2$ or if $\lambda = 2 \cos(\pi/m)$,
where $m \geq 3$ is an integer, $B(\lambda)$ is
a fundamental region for $G(\lambda)$.
\end{theorem}
\begin{definition}
\footnote{\cite{berndt2008hecke}, Definition 3.4} 
Let $T_A = \{ \lambda: 
\lambda = 2 \cos(\pi/m), m \geq 3, m \in \mathbb{Z}\}$.
\end{definition}
\noindent
Berndt states in his Theorem 5.4 that 
$G(\lambda)$ is discrete if and only if
$\lambda$ belongs to $T_A$. This discreteness
is the premise of the theory of 
automorphic functions generally.
He embeds within the
proof of his Lemma 3.1 
(which we omit), the 
\begin{definition}
$\tau_{\lambda}$ denotes the intersection
in $\mathbb{H}$ of the line
$x = -\lambda/2$ and the unit circle
$|\tau| = 1$. 
\end{definition} \noindent
(Berndt remarks at the top of page 35 that
$\tau_{\lambda}$ is the lower left corner
of $B(\lambda)$). and that
$\pi \theta = \pi - \arg(\tau_{\lambda})$,
so that $\cos (\pi \theta) = \lambda/2$.)
\newline \newline \noindent
To characterize Eisenstein series,
we need to keep track of some analytical
properties. The next definition summarizes
the second paragraph  of
Berndt's Chapter 5.
(Throughout his Chapter 5, $\lambda < 2$.)
\begin{definition}
Let $f \in M(\lambda, k, \gamma), f$
not identically zero.
    \begin{enumerate}
    \item
    $N = N_f$ counts the zeros of $f$ on 
    $\overline{B(\lambda)}$  with multiplicities.
    \item $N_f$ does not count zeros
    at $\tau_{\lambda}$, at $\tau_{\lambda} + \lambda$, 
    at $i$, or at $i\infty$.
    \item 
    If $\tau_0 \in 
    \overline{B(\lambda)}, f(\tau_0) = 0$
    and
   $\Re (\tau_0) = -\lambda/2$, 
    then $f(\tau_0 + \lambda) = 0$ and
    $N_f$ counts only one of the two zeros.
    \item If $\tau_0 \in 
    \overline{B(\lambda)}, f(\tau_0) = 0$,
    and $|\tau_0| = 1$, then,
     $f(-1/\tau_0) = 0$,
    and $N_f$ counts only one 
    of these two zeros.
    \item  The numbers $n_{\lambda}, n_i,$
    and $n_{\infty}$ are the orders of the
    zeros of $f$ at $\tau_{\lambda},
    i$ and $i\infty$, repectively.
    The order $n_{\infty}$ is measured in terms of
    $\exp(2 \pi i \tau/\lambda)$.
\end{enumerate}
\end{definition}
\noindent
The multiplier $\gamma$ is given by
\begin{theorem}
\footnote{\cite{berndt2008hecke}, Corollary 5.2} 
Let $f \in M(\lambda,k,\gamma)$ and 
let $n_i$ be the order of the zero
of $f$ at $\tau = i$.
Then
$$
\gamma= (-1)^{n_i}.
$$
\end{theorem}
\noindent
The next two results 
tell us that the only nontrivial
case in this theory is the one 
that we are interested in.
\begin{theorem} 
\footnote{\cite{berndt2008hecke}, Lemma 5.1} 
If $\dim M(\lambda,k,\gamma)
\neq 0$,
$$
N_f + n_{\infty} + \frac 12 n_i + 
\frac {n_{\lambda}}m =
\frac 12 k \left ( \frac 12 - \theta \right).
$$
\end{theorem}
\noindent
By Berndt's equation (5.16), 
if $m \geq 3$ then 
the right side can be written as $k(m-2)/4m$.
\begin{theorem}
\footnote{\cite{berndt2008hecke}, Theorem 5.2}
If $\dim M(\lambda,k,\gamma)
\neq 0$, then $\theta = 1/m$
where $m \geq 3$ and $m \in \mathbb{Z}$.
\end{theorem}
\noindent
We are concerned with
 $\lambda \in T_A$. 
This makes $\lambda < 2$ as in all 
the results of Berndt's 
Chapter 5.
\newline \newline \noindent
One estimate 
for $\dim M(\lambda, k, \gamma)$ is 
\begin{theorem}
\footnote{\cite{berndt2008hecke}, Theorem 5.6}
If $\lambda \notin T_A$, then
$\dim M(\lambda,k,\gamma) = 0$.
If $\lambda = 2\cos(\pi/m) \in T_A$,
then for nontrivial $f \in
M(\lambda,k,\gamma)$, the weight
$k$ has the form
$$
k = \frac {4h}{m-2} + 1 - \gamma,
$$
where $h \geq 1$ is an integer.
Furthermore,
$$
\dim M(\lambda,k,\gamma) = 1 + \left \lfloor
\frac{h + (\gamma-1)/2}m \right \rfloor.
$$
\end{theorem}
\noindent
Eliminating $h$, we find that
\begin{equation}
\dim M(\lambda,k,\gamma) = 1 +
\left \lfloor 
k\left (\frac 14 - \frac 1{2m} \right ) + 
\frac {\gamma}4 - 
\frac 14
\right \rfloor.
\end{equation}
Berndt
proves that the  dimension formula 
 above holds also when $h = 0$.
 \footnote{\cite{berndt2008hecke}, Remark 5.3}  
 \newline  \newline \noindent
The existence of 
 certain modular forms is
 provided by
 \begin{theorem}
\footnote{\cite{berndt2008hecke}, Theorem 5.5}
Let $\lambda \in T_A$. Then there exist
 functions $f_{\lambda}, f_i$, and 
 $f_{\infty} \in M(\lambda,k,\gamma)$
 such that each has a simple zero at 
 $\tau_{\lambda}, i$, and $i \infty$,
 respectively, and no other zeros. 
 Here, $\gamma$ is given by 
 Theorem 3 of the present article,
 and $k$ is determined in each case from
 Theorem 4 of the present article.
 Thus, $f_{\lambda} \in 
 M(\lambda, 4/(m-2), 1), f_i \in 
 M(\lambda, 2m/(m-2), -1)$, and
 $f_{\infty} \in M(\lambda, 4m/(m-2),1)$.
 \end{theorem}
 \begin{remark}
 \footnote{\cite{berndt2008hecke}, pages 47-48}
By the Riemann mapping theorem
there exists a function $g(\tau)$ that
maps the simply connected region $B(\lambda)$
one-to-one and conformally
onto $\mathbb{H}$. If
we require that
$g(\tau_{\lambda}) = 0,
g(i) = 1$, and $g(i\infty) = 
\infty$, then $g$ is determined uniquely.
 \end{remark}
 \noindent
 Now we can write down 
 $f_{\lambda}, f_i$, and $f_{\infty}$ explicitly.
 The next theorem is extracted from 
the proof of Theorem 7.
$f_{\lambda}$ and $f_i$
correspond to Eisenstein series and 
$f_{\infty}$ to a cusp form.
In our code, we take
 $g$ to be a normalized form of $J_m$.
 \begin{theorem}\footnote{\cite{berndt2008hecke}, page 50}
 $$
 f_{\lambda}(\tau)
 =\left \{
 \frac {g'(\tau)^2}
 {g(\tau)(g(\tau) -1)}
 \right \}^{1/(m-2)},
 $$
 $$
 f_i(\tau) =
 \left \{
 \frac {g'(\tau)^m}
 {g(\tau)^{m-1} (g(\tau) - 1)}
 \right \}^{1/(m-2)},
 $$
 and
 $$
 f_{\infty}(\tau) = 
 \left \{ \frac{g'(\tau)^{2m}}
 {g(\tau)^{2m-2}(g(\tau)-1)^m}
\right \}^{1/(m-2)}.
 $$
 \end{theorem}
 \noindent
In our applications to Lehmer's
problem, we will be
interested in the dimensions
of the weight $12$ cusp spaces
for $\lambda = \lambda_m = 2 \cos \pi/m$.
\begin{definition}
\footnote{\cite{berndt2008hecke}, Definition 5.2}
If
$f\in M(\lambda, k, \gamma)$
and $f(i \infty) = 0$,
then we call $f$ a cusp form of
weight $k$ and multiplier  $\gamma$
with respect to $G(\lambda)$. 
We denote by $C(\lambda,k,  \gamma)$ the
vector space of  all cusp forms
of this kind.
\end{definition}
\begin{remark}
\footnote{\cite{berndt2008hecke}, equation (5.25)}
$$\dim C(\lambda, k,  \gamma) \geq 
\dim M(\lambda, k,  \gamma) -1.$$
\end{remark}
\begin{remark}
In view of (i) Theorem 6, (ii) equation (12), 
(iii) Remark 2, and (iv)  the fact
that $\gamma = \pm 1$, we see that
$\dim C(\lambda_m, 12,  \gamma) >1$
when $m$ is greater than or equal to $12$.
\end{remark}
\section{Modular forms studied 
in our experiments}
We are going to
write down versions of 
the functions from Theorem 8
such that, 
at $m = 3$, they 
reduce  to corresponding
functions in the classical theory.
Some have fixed weights (four,
six and twelve)
and others have weights that vary with $m$.
The classical objects
(in Serre's notation \cite{serre1970course}),
are Klein's $j$-invariant, the weight four
Eisenstein series $E_2$,
the weight six Eisenstein series $E_3$,
and the generating function of Ramanujan's
tau function, namely the normalized weight twelve
cusp form $\Delta$. 
They all belong to one-dimensional
vector spaces of modular forms
and the number of zeros each one has in a 
given fundamental region is small,
so  the identifications
follow by comparison of the
initial Fourier coefficients.
\footnote{\cite{serre1970course}, Chapter VII, equations (20-21)}
\newline \newline \noindent
Corresponding to $f_{\lambda}$,
we have
\begin{definition}
\begin{enumerate}
\item
$H_{\lambda, m}(\tau)$:=
 $$\left \{
 \frac {J_m'(\tau)^2}
 {J_m(\tau)(J_m(\tau) -1)}
 \right \}^{1/(m-2)}.$$
 \item $H_{\lambda,4,m} (\tau) := 
 H_{\lambda, m}(\tau)^{m-2}.$
 \end{enumerate}
 \end{definition}
 \noindent
  Corresponding to $f_i$, 
 we state
 \begin{definition}
 \begin{enumerate}
 \item
  $$
 H_{i,m}(\tau):=
\left \{
 \frac {J_m'(\tau)^m}
 {J_m(\tau)^{m-1} (J_m(\tau) - 1)}
 \right \}^{1/(m-2)}.$$
 \end{enumerate}
 \end{definition} 
 \begin{definition}
 \begin{enumerate}
 \item
  Corresponding to $f_{\infty}$,
 we have
 $$
 \Delta_{\infty,m}(\tau) := 
 \left \{ \frac{J_m'(\tau)^{2m}}
 {J_m(\tau)^{2m-2}(J_m(\tau)-1)^m}
\right \}^{1/(m-2)}.
 $$
 \item $\Delta^{\diamond}_m :=
H_{\lambda,m}^3/J_m$.
\item 
$\Delta^{\diamond}_{12,m} :=
H_{\lambda,4,m}^3/J_m$.
 \end{enumerate}
 \end{definition} 
 \begin{remark}
  By Berndt's theorem 7 above,
  we have the following table of weights:
 \begin{center}
\begin{tabular}{|c|c|c|c|c|c|} \hline
$H_{\lambda,m}$& $H_{\lambda,4,m}$ &$H_{i,m}$
& $\Delta^{\diamond}_m$ 
&  $\Delta^{\diamond}_{12,m}$ & $\Delta_{\infty,m}$ \\ \hline
 $4/(m-2)$ &  $4$  & $2m/(m-2)$ &  $12/(m-2)$ &  $12$ &
 $4m/(m-2)$ \\ \hline
\end{tabular} 
\end{center}
 \end{remark}
\section[]{Interpolation by polynomials}
In this section, we state conjectures about
polynomials interpolating coefficients of
modular forms for Hecke groups.
Conjectures 6 and 7
bear on Lehmer's question about
the existence of zeros of Ramanujan's
tau function.
\newline \newline \noindent
Berndt's (Hecke's) theorems
7 and 8 above make it clear that 
 Akiyama's theorem proving
 Raleigh's conjecture
 on the interpolation of
 the coefficients of the
 Fourier expansions of 
 Hecke triangle functions
extends in some way to the 
modular forms defined in
the previous section. We
did experiments  to
explore the details; 
our observations are summarized
in the conjectures below.
\subsection[]{Analogues of $SL(2,\mathbb{Z})$ 
Eisenstein series.}
We found the sequence
$\{e_{4,n}\}$ mentioned below on \cite{OEIS14}.
\begin{conjecture}
\footnote{Notebooks ``conjecture 3.nb''
and ``conjecture 3.ipynb'',
\cite{test3};
associated data files in the data folder on
\cite{test3}.}
Let the Fourier expansion of 
$\overline{H_{\lambda,4,m}}(\tau)$ be 
$$\overline{H_{\lambda,4,m}}(\tau)=
\sum_{n=0}^{\infty} \beta_{4,m}(n) X_m^n.$$
\begin{enumerate}
\item \footnote{\cite{serre1970course}, 
page 93.} $\overline{H_{\lambda,4,3}}(\tau)$
reduces to Serre's weight-4 Eisenstein series
$E_2$ in the sense that 
$\beta_{4,3}(n) = 240 \sigma_3(n)$ for 
$n = 1, 2, 3, ...$.
\item For each $n$ there is a polynomial
$B_{4,n}(x)$ with rational coefficients
such that 
$m^{3n}\beta_{4,m}(n)=B_{4,n}(m)$ for  $m = 3,4, ....$
\item If $n$ is positive, then the degree
of $B_{4,n}(x)$ is $6n$.
\item $B_{4,0}(x) \equiv 1$ and,
if $n$ is positive, then 
$$
B_{4,n}(x) = 
e_{4,n}  (x^2-4)x^{4n}  b_{4,n}(x),
$$
where $e_{4,n} = 
16 \sum_{\substack{\nu | n \\ \nu \text{odd}}} 
(-1)^{n-\nu} \nu^3$
and $b_{4,n}(x)$ is a monic irreducible polynomial
in $\mathbb{Q}[x]$.
\end{enumerate}
\end{conjecture}
\begin{conjecture}
\footnote{Notebooks ``conjecture 4.1-4.3.ipynb,
conjecture 4.4a.ipynb,
conjecture 4.4b.ipynb'',
conjecture 4.5.ipynb,
\cite{test3}. N.B.:
Contrary to appearances,
the function denoted ``H4''
in these \it SageMath \rm
notebooks is not
the function covered in the
previous conjecture.
``H4'' is $H_{\lambda,m}$.
}
Let the Fourier expansion of 
$\overline{H_{\lambda,m}}$ be
$$
\overline{H_{\lambda,m}} =
\sum_{n=0}^{\infty} \beta_m(n) X_m^n.
$$
\begin{enumerate}
\item For each $n$ there is a polynomial
$B_n(x)$ with rational coefficients
such that 
$\beta_m(n)=B_n(m)$ for  $m = 3,4, ....$
\item If $n$ is positive, then the degree
of $B_n(x)$ is $3n - 1$.
\item $B_0(x) \equiv 1$ and $B_1(x)  = 16x(x+2)$.
\item Let  $\mathscr{Q}$ be as in item 9 of
our glossary and let 
$e_n =  16 (-1)^{n+1} 
\sum_{\substack{\nu | n \\ \nu \text{odd}}} 1/\nu$.
If $n$ is greater than $2$ and belongs to 
$\mathscr{Q}$, then 
$$
B_n(x) = 
e_n  (x^2-4)(x-6)x^n  b_n(x),
$$
where $b_n(x)$ is a monic irreducible polynomial.
Otherwise (for $n$ greater than one)
$B_n(x) =  e_n (x^2-4)x^n  b_n(x)$
where, again,
$b_n(x)$ is a monic irreducible 
polynomial in $\mathbb{Q}[x]$.
\item $\overline{H_{\lambda,3}}$
reduces to $E_2$ in the same sense as
in conjecture 3.1.
\end{enumerate}
\end{conjecture}
\noindent
(We identified the $e_n$ after reading
\cite{OEIS11} and \cite{OEIS12}.)
\newline \newline \noindent
Thus, in the range of our observations
($3\leq m\leq 302, 0\leq n\leq 100)$,
the only integer value of $m$
such that
$\overline{H_{\lambda,m}}$ has any vanishing coefficients
is six, and $\beta_n(6)$ is zero just if
$n$ is in $\mathscr{Q}$. 
\begin{conjecture}
\footnote{In \cite{test3} notebook ``conjecture 5.ipynb''.}
Let the Fourier expansion of $\overline{H_{i,m}}$ be
$$
\overline{H_{i,m}} =
\sum_{n=0}^{\infty} \delta_m(n) X_m^n.
$$ 
\begin{enumerate}
\item For each non-negative integer $n$, 
there is a polynomial
$D_n(x)$ in $\mathbb{Q}[x]$ 
such that 
\begin{enumerate}
    \item $D_n(m) = \delta_m(n)$ for $n = 0, 1,...$
    and $m  = 3, 4,....$
    \item The degree of $D_n$ is $3n$.
    \item $D_n(x) =$ a rational number 
    $d_n \times$ a product of 
    monic irreducible polynomials.
    \item $d_0 = 1$  and, for $n$ a positive
    integer, $d_n = 24(-1)^n  \sum^*_{\nu|n} \nu$.
    Again, the asterisk means that the
    sum is taken over the odd positive divisors of $n$.
\end{enumerate} 
\item  $D_n(m)=(-1)^m \delta_n(m)$  for  $m = 3,4, ....$
\item $D_0(x) \equiv 1$ identically,
$D_1(x)  = -24(x-2/3)x^2$, and
$D_2(x) = $ \newline $24(x-2/3)(x-2)x^3(x-14)$.
\item For $n$ larger than two,
$D_n(x) =d_n (x-2)(x-2/3)x^{n+1} \epsilon_n(x)$
where $\epsilon_n(x)$ is a monic irreducible
polynomial in  $\mathbb{Q}[x]$.
\item $\overline{H_{i,3}}$ reduces to Serre's weight-6
Eisenstein series $E_3$ in the sense that
$\delta_3(0) = 1$ and $\delta_3(n) = -504\sigma_5(n)$
for $n=1,2,3,...$.
\end{enumerate}
\end{conjecture}
\subsection[]{Analogues of $SL(2,\mathbb{Z})$ cusp forms.}
Let $\Delta$ be usual normalized discriminant,
 a weight $12$ cusp form
for $SL(2,\mathbb{Z}) =  G(\lambda_3)$
with integer coefficients.
Its Fourier expansion is written
$$
\Delta(\tau) = \sum_{n=1}^{\infty} \tau(n) q^n
$$
where $q = e^{2 \pi i \tau}$ and
 $\tau(n)$ is Ramanujan's function. 
 (The reader will not 
 confuse the complex number $\tau$
 with Ramanujan's function $\tau(n)$
 or any of its relatives defined below.)
Whether or not
the equation $\tau(n) = 0$ has 
any solutions is,
of course, an open
question.\footnote{\cite{lehmer1947vanishing}}
Several authors have eliminated
various classes of integers as values of 
tau.\footnote{These results 
are summarized in 
\cite{lakein2021some}.
Relevant citations are 
\cite{balakrishnan2020variations},
\cite{balakrishnan2021even},
\cite{balakrishnan2020variants},
\cite{bennett2021odd},
\cite{dembner2021hyperelliptic},
\cite{lygeros2013odd},
\cite{murty1987odd},
and \cite{lakein2021some}
itself.
} 
It  will be apparent that each of the conjectures
about cusp-form analogues implies that tau has no zeros.
\newline \newline \noindent
From definition 10.2,
 $$
 \Delta_{\infty,m}(\tau)^{m-2} = 
 \frac{J_m'(\tau)^{2m}}
 {J_m(\tau)^{2m-2}(J_m(\tau)-1)^m}
 $$
 and, by Theorem 7 in our sketch of
 Hecke's theory, its weight is $4m$. Since it
 raises a cusp form beginning with an $X^1$ term
 to high powers, we will use the star operator
 (definition 1.2) to state the following
 conjecture.
\begin{conjecture}\footnote{Notebooks 
``conjecture 6Laptop.nb'',
``conjecture 6.ipynb''
and ``conjecture 6 no2.ipynb'',
\cite{test3}.}
 Let the Fourier expansion
 of $\overline{(\Delta_{\infty,m}(\tau)^{m-2})^*}$ 
 be written as
 $\overline{(\Delta_{\infty,m}(\tau)^{m-2})^*} = 
 \sum_{n=0}^{\infty} \overline{\tau}_m(n) X_m$.
\begin{enumerate}
    \item $\overline{\tau}_3(n-1) = \tau(n)$ for $n = 1, 2, ...$.
     \item There is a set of polynomials
     $\overline{T}_n(x), n = 1, 2, 3...$
     such that, for each $n, 
     \overline{T}_n(m) =\overline{\tau}_m(n)$.
     \item
     $\overline{T}_n(x) = (-8)^n (x-2)^3 x^n t_n(x)/n!$
     where $t_n$ is a polynomial
     with rational coefficients that is irreducible 
     over $\mathbb{Q}[x]$.
 \end{enumerate}
\end{conjecture}
\begin{conjecture}\footnote{Notebook 
``conjecture 7.ipnyb'',
\cite{test3}}
Let the Fourier expansion of 
$\overline{\Delta_{\infty,m}}(\tau)$ be
$$
\overline{\Delta_{\infty,m}}(\tau) =
\sum_{n=1}^{\infty} \tau_{\infty,m}(n) X_m^n.
$$
\begin{enumerate}
\item $\tau_{\infty,3}(n) = \tau(n)$
for $n = 1,2, 3,...$.
\item There is a set of polynomials
$T_{\infty,n}(x)$
with coefficients in $\mathbb{Q}$ such that 
$\tau_{\infty,m}(n) = T_{\infty,n}(m)$.
\item
$T_{\infty,1}(x) \equiv 1$ identically,
and, if $n$ is greater than one,  
\begin{enumerate}
\item $T_{\infty,n}(x) = 
s_{\infty,n}(x-2)^2 x^{n-1} t_{\infty,n}(x)$,
where $t_{\infty,n}(x)$ is a
monic irreducible polynomial 
over $\mathbb{Q}$
of degree $2n - 4$ and
\item
$s_{\infty,n}$ is 
(in the notation of \cite{conway2013sphere},
 Chapter 7, Theorem 7) the
coefficient of $q^n$ in the Fourier
expansion of $\Delta_8(z)$.\footnote{\cite{OEIS3}}
\item \footnote{\cite{OEIS4}}
Also, 
$$
s_{\infty,n} = (-1)^{n+1} 
\sum_{\substack{\nu | n \\ n/\nu \text{odd}}} \nu^3.
$$
This sum is the coefficient of 
$q^n$ in the Fourier expansion
of $E_{\infty, 4}$,
the unique normalized weight-$4$ modular form 
for $\Gamma_0(2)$ with simple zeros at 
$i \infty$);\footnote{\cite{brent1998quadratic},
equation (2-3)}
it is also the number of representations of 
$n-1$ as a sum of $8$ triangular numbers.
\footnote{\cite{ono1995representation}, Theorem 5.}
\item
Finally, $s_{\infty,n}$ is 
the coefficient of $q^n$ in the expansion 
of $\eta(2z)^{16}/\eta(z)^{-8}$
where $\eta(z)$ is Dedekind's function
(\cite{brent1998quadratic}, equation (2-16).)
\end{enumerate}
\end{enumerate}
\end{conjecture}
\begin{conjecture}\footnote{Notebook 
``conjecture 8.ipynb'', \cite{test3}.}
Let the Fourier
expansion of 
$\overline{\Delta^{\diamond}_m}(\tau)$ be
$$
\overline{\Delta^{\diamond}_m} = 
\sum_{n=1}^{\infty}  
\tau_m^{\diamond}(n) X_m^n.
$$
\begin{enumerate}
\item $\tau_3^{\diamond}(n) = \tau(n)$
for $n = 1,2,3,...$.
\item There is a set of polynomials
$T^{\diamond}_n(x)$
with coefficients in $\mathbb{Q}$ such that 
$\tau_m^{\diamond}(n) = 
T^{\diamond}_n(m)$.
\item $T^{\diamond}_1(x), 
T^{\diamond}_2(x)$, and 
$T^{\diamond}_3(x)$
are irreducible polynomials over 
$\mathbb{Q}$ of degrees
$3, 6$, and $9$, respectively.
\item If $n$ is greater than $3$,
$T^{\diamond}_1(x) = 
s^{\diamond}_n \cdot (x-2)x^{n-1} 
t^{\diamond}_n(x)$, 
where
$s^{\diamond}_n$ is a rational number and 
$t^{\diamond}_n(x)$ is a monic
polynomial, irreducible 
over $\mathbb{Q}$, of degree $2n-3$.
Furthermore,
\begin{enumerate}
\item
$\sum_{n=0}^{\infty}
s^{\diamond}_n q(\tau)^n
=  $ 
$$
 \prod_{n \text{odd}}(1-q(\tau)^n)^{24} \times
\prod_{n \equiv 2 (4)}(1-q(\tau)^n)^{-24} = 
\eta^{24}(\tau) \eta^{24}(4\tau) \eta^{-48}(2\tau).
$$
\item
$s^{\diamond}_n = 
(-1)^{n+1} \times$ the coefficient 
of $q(\tau)^n$ in $(\eta(2\tau)/\eta(\tau))^{24}$.
\end{enumerate}
\item There is no corresponding set of
interpolating polynomials for 
$\Delta^{\diamond}_3$.
\end{enumerate}
\end{conjecture}
\noindent
The product decomposition in clause 3(a) above
is a guess based on $43$ terms of the series
using Euler's method.
\footnote{\cite{apostolIntroduction} (Theorem 14.8);
 English-language version of \cite{eulerdiscovery}
 in \cite{polya2021mathematics};
and \cite{OEIS6}. The second decomposition
appears in \cite{OEIS7}.}
\begin{conjecture}
 \item \footnote{Notebook ``conjecture 9.ipynb'',
\cite{test3}.}
Let the Fourier expansion of 
$\overline{\Delta^{\diamond}_{12,m}}(\tau)$ be
$$
\overline{\Delta^{\diamond}_{12,m}}(\tau) = 
\sum_{n=1}^{\infty} 
\tau_{12,m}^{\diamond}(n) X_m^n.
$$
\begin{enumerate}
\item 
$\tau_{12,3}^{\diamond}(n) = \tau(n)$ for 
$n = 1, 2, 3 ...$.
\item There is a set of polynomials
$T^{\diamond}_{12,n}(x), n = 1, 2, ...$ of degree $3n-3$
with coefficients in $\mathbb{Q}$ such that 
$\tau_{12,m}^{\diamond}(n) = 
T^{\diamond}_{12,n}(m)$
for each $m = 3, 4, ....$
\item For each $n$, there are zeros of 
$T^{\diamond}_{12,n}(x)$
on both axes of the complex plane,
and there are no other complex
zeros.\footnote{Notebook
``conjecture 9.nb'' \cite{test3}
contains plots of the complex zeros for
$n$ between $1$ and $24$.}
(Figures 2 and 3 illustrate this for
$n = 11$ and $24$.)
\item 
$T^{\diamond}_{12,n}(x) = 
(-1)^{n+1}\tau(n) x^{n-1} t^{\diamond}_{12,n}(x)$,
where  $t^{\diamond}_{12,n}(x)$ 
is a monic irreducible polynomial
over $\mathbb{Q}$.
\end{enumerate}
\end{conjecture}
\section[]{Lehmer's question}
 \begin{remark} By clause 4 of conjecture 9, 
 for $m = 3, 4, ...$\thinspace:
 $\tau^{\diamond}_{12,m}(n) = 0$
 if and only if $\tau(n) = 0$. 
 \end{remark} \noindent
More generally, we have
\begin{conjecture}
Letting $T_n(x)$ and  $\tau_m$ stand for
the various polynomials and Fourier coefficients
in conjectures 6 through 9,
none of the $T_n(x)$ has an integer 
root greater than two; 
consequently, none of the $\tau_m$ vanishes
for $m = 3, 4, ....$
\end{conjecture} \noindent
 Let $d(m,n)$ be the minimum 
 Euclidean distance to $m$
 of any complex root of $T_n(x)$.
 We have (in effect)
 conjectured above that
 in each case $T_n(3) = \tau(n)$,
 so the behavior of $d(3,n)$
 measures how closely we
 can come to the assertion
 that $\tau(n) = 0$ for some $n$.
\begin{conjecture}
 For any positive real number $r$,
  $d(3,n)$ is less that $e^{-rn}$
 for sufficiently large $n$.\footnote{For 
 this proposal, we depend on  
 graphical evidence
 which we sample figures 10 -- 17. More
 extensive collections of plots are
 in notebooks ``conjecture 6.1.nb'',
 ``conjecture 6.2.nb'', ``conjecture 7.nb'',
 and ``conjecture 8.nb'', \cite{test3}.} )
\end{conjecture}
\section{Other questions}
\begin{enumerate}
\item 
Like $G_n$ in clause 5 of conjecture 1, 
the index-$n$ hyperoctahedral group
has size  $2^n n!$\thinspace.
\footnote{\cite{hyper,miller1918groups,graczyk2005hyper}}
Are they isomorphic?
\item In conjectures 1--9,
the $n^{th}$ interpolating polynomial
is written as a product of a numerical term
and several monic polynomials belonging to 
$\mathbb{Q}[x]$. In each case,
all but one of the monic factors is 
given explicitly, \it i.e.\rm,
in terms of $n$, but without reference
to the Fourier expansion of
the underlying modular form. The
``inexplicit'' factor
can, of course, be written in terms of
the first $n$ of these coefficients, 
but can it be expressed in the same way as 
the other factors: without reference to the 
Fourier coefficients?
\item
While checking our calculations, we compared the 
Fourier expansion
of $H_{\lambda,4}(x/A_4)$ (abusing notation in the 
obvious way) with Leo's 
expansion of the weight $4$ Eisenstein series at $m = 4$.
\footnote{\cite{leo2008fourier}, p.54}
(Recall that $A_4 = 1/256.$)
Within the range of our observations, they do coincide. 
%
%
The expansions (in our own notation) both begin 
$$
1 + 48q_{_4} + 624 q_{_4}^2 + 1344 q_{_4}^3 + ....
$$
Let 
$$
E_{\gamma,2} = 
\newline
1+24\sum_{n=1}^{\infty} 
\sum_{\substack{\nu | n \\ \nu \text{odd}}} \nu q^n.
$$
Sloane comments that the sequence
$\{1,48, 624, ... \}$ is the same as that of the
coefficients of $E_{\gamma,2}^2$.\footnote{\cite{OEIS5}}
$E_{\gamma,2}^2$ is a weight $4$, 
level $2$ modular form, that is,
a weight $4$ modular form for the 
$SL(2,\mathbb{Z})$ subgroup
$\Gamma_0(2)$.\footnote{\cite{brent2001quadratic},
equation (2-1), p, 260}.
We propose in conjecture 7 (c) above
that $s_{\infty,n}$ is the coefficient of 
$q^n$ in the Fourier expansion
of $E_{\infty, 4}$,
the unique normalized weight-$4$ modular form 
for $\Gamma_0(2)$ with simple zeros at 
$i \infty$.
We have also proposed in conjectures 1, 2,
7 and 8 that
interpolating polynomials 
are products of monic polynomials
with rational numbers
equal or related to
Fourier coefficients of
other classical Hauptmoduln. 
What  is the relationship
between modular forms for 
subgroups of $SL(2,\mathbb{Z})$ and  
modular forms for the
other $G(\lambda_m)$?
\item Both $J_m$ and $\overline{J_m}$
(that is, $j_m$) appear to be interpolated
by polynomials. On the other hand,
$\overline{\Delta^{\diamond}_m}$ 
appears to be interpolated by 
polynomials, but 
$\Delta^{\diamond}_m$ does not.
Why are the situations different?
\end{enumerate}
\section{Figures}
\subsection{Figure 1.} 
\hskip 0in
\includegraphics[scale=.7]{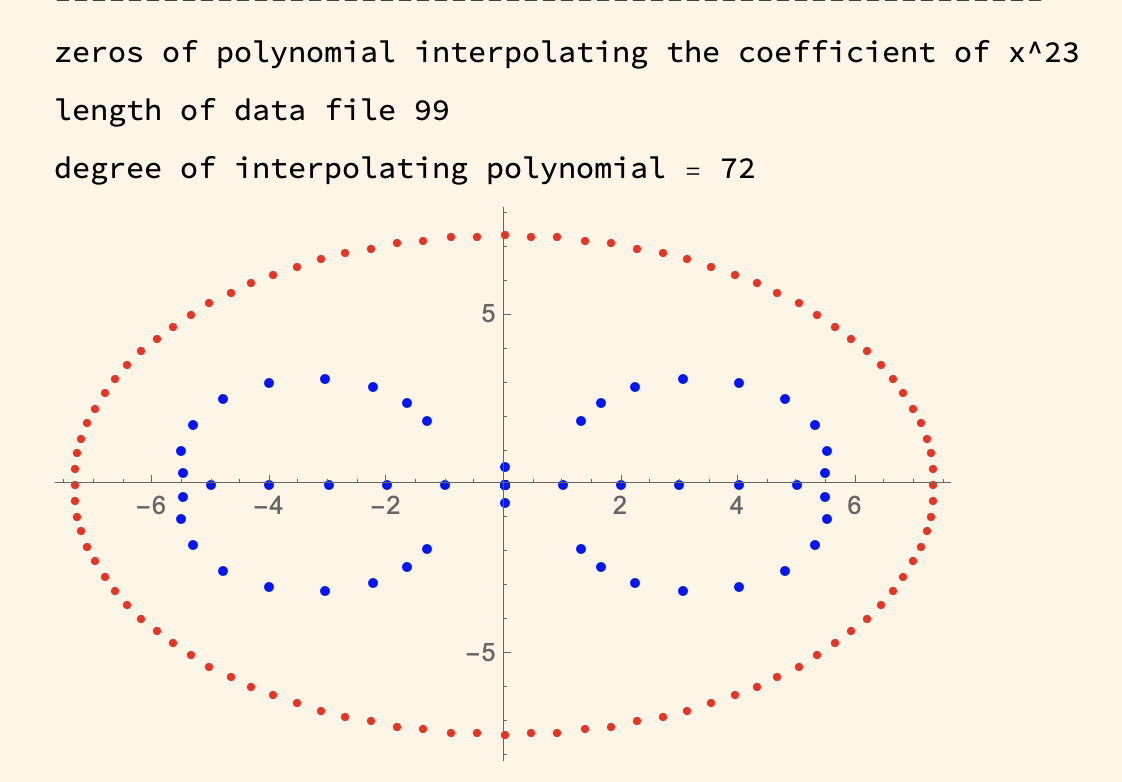} 
\sc{Roots of polynomial interpolating 
the coefficient of $q_m^{23}$
in the Fourier expansion of $j_m(\tau)$
 (conjecture 1.)
 \footnote{Notebook 
``conjecture1clause4d.nb'',
\cite{test3}.}}
\subsection{Figure 2.}
\includegraphics[scale=.8]{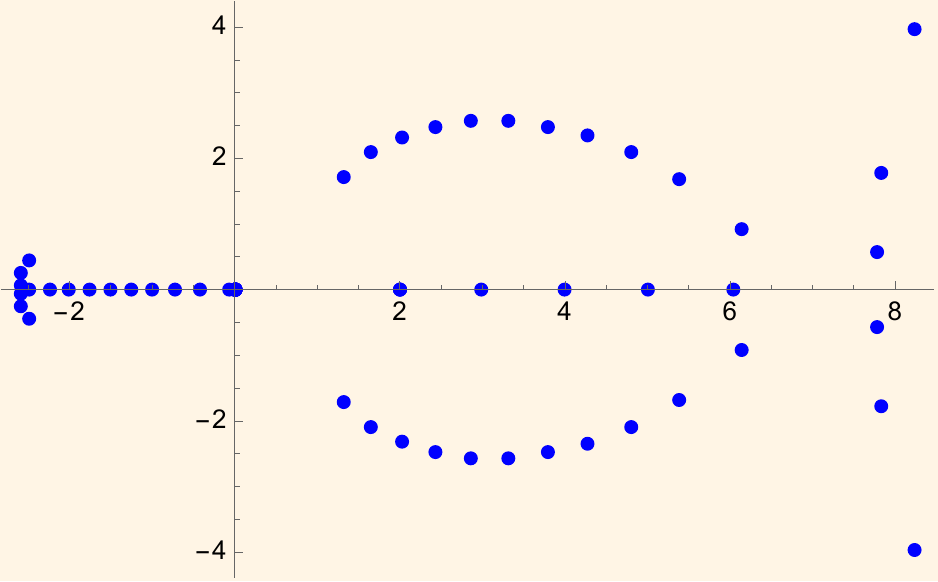}
\sc{Roots of $\overline{T}_{17}$ (conjecture 6.)
\footnote{Notebook ``conjecture 6Laptop.nb'',
\cite{test3}.}}
\subsection{Figure 3.}
\includegraphics[scale=.8]{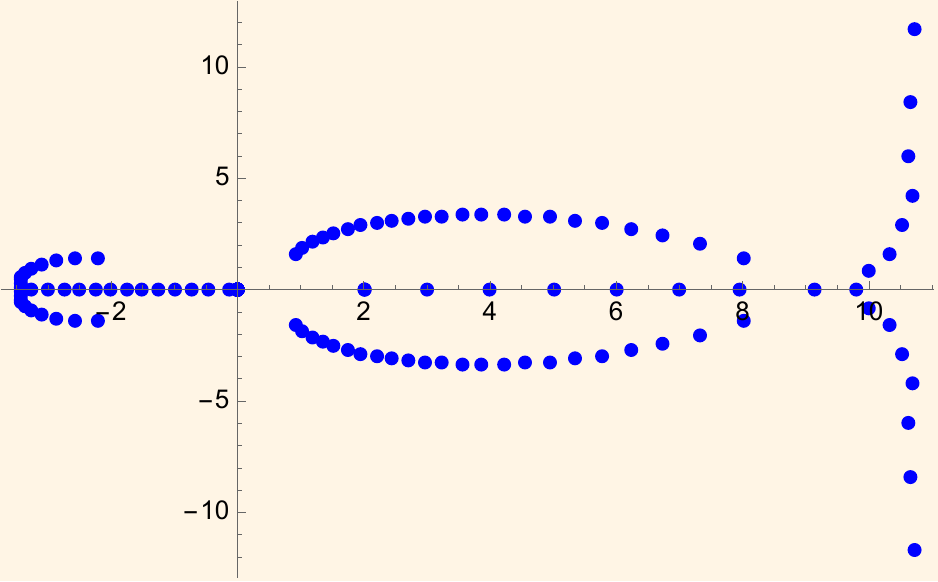}
\sc{Roots of $\overline{T}_{35}$ (conjecture 6)
\footnote{\it ibid.\rm}}
\subsection{Figure 4.}
\includegraphics[scale=1]{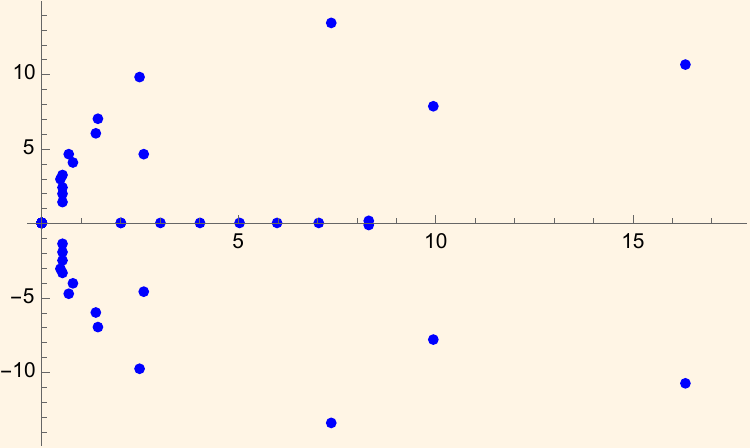}
\sc{Roots of $T_{\infty,20}$ (conjecture 7.)
\footnote{Notebook ``conjecture 7.nb'',
\cite{test3}}}
\subsection{Figure 5.}
\includegraphics[scale=1]{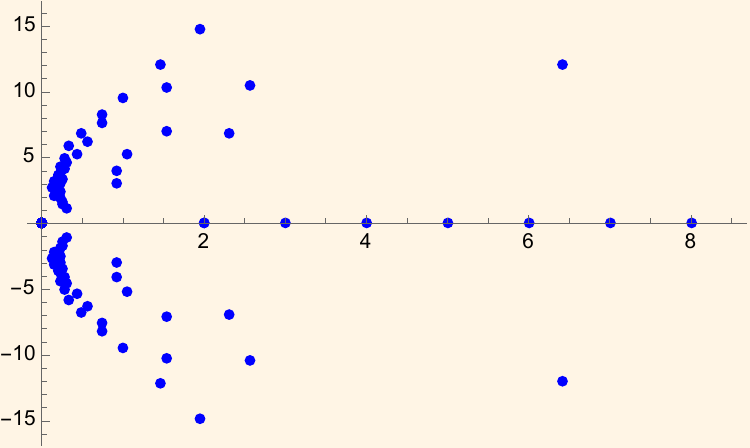}
\sc{Roots of $T_{\infty,50}$ (conjecture 7.)
\footnote{\it ibid.\rm}}
\subsection{Figure 6.}
\includegraphics[scale=1]{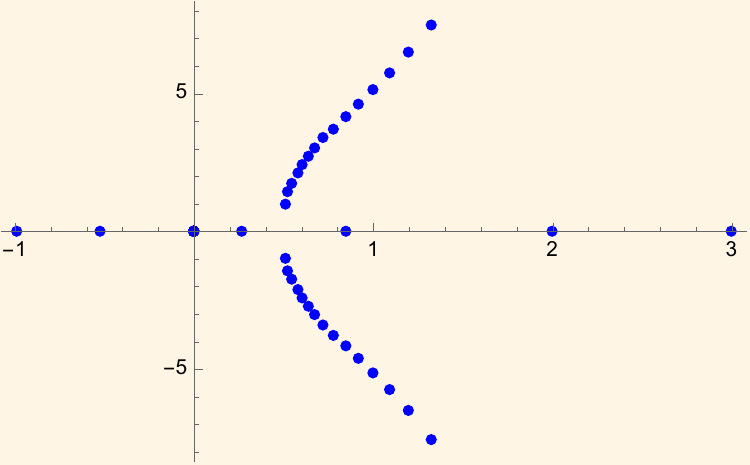}
\sc{Roots of $T^{\diamond}_{19}$ (conjecture 8.)
\footnote{Notebook ``conjecture 8.nb'',\cite{test3}.}}
\subsection{Figure 7.}
\includegraphics[scale=1]{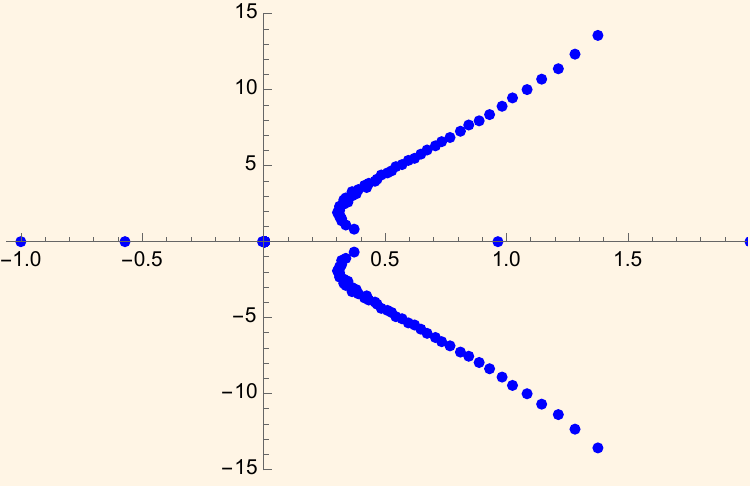}
\sc{Roots of $T^{\diamond}_{50}$ (conjecture 8.)
\footnote{\it ibid.\rm}}
\subsection{Figure 8.}
\includegraphics[scale=1]{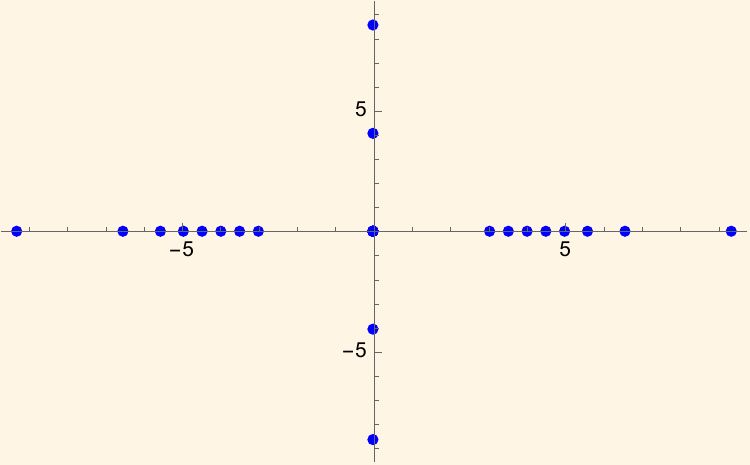} 
\sc{Roots of $T^{\diamond}_{12,11}$ (conjecture 9.)
\footnote{Notebook ``Conjecture 9.nb'', 
\cite{test3}.}}
\subsection{Figure 9.}
\includegraphics[scale=1]{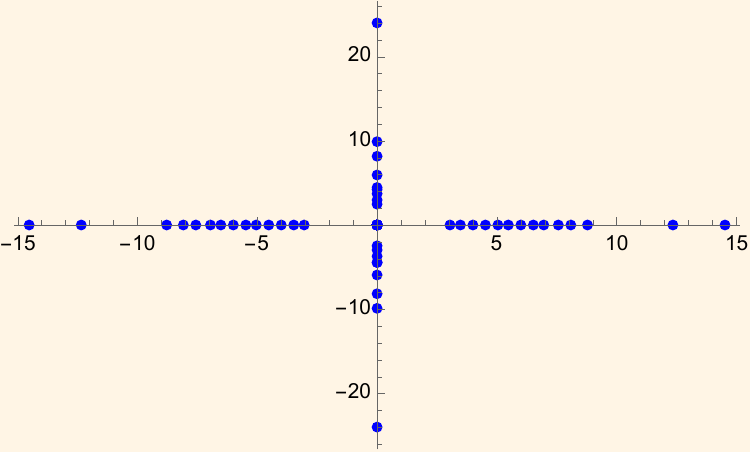}
\sc{Roots of $T^{\diamond}_{12,24}$ (conjecture 9.)
\footnote{\it ibid.\rm}}
\subsection{Figure 10.}
\includegraphics[scale=1]{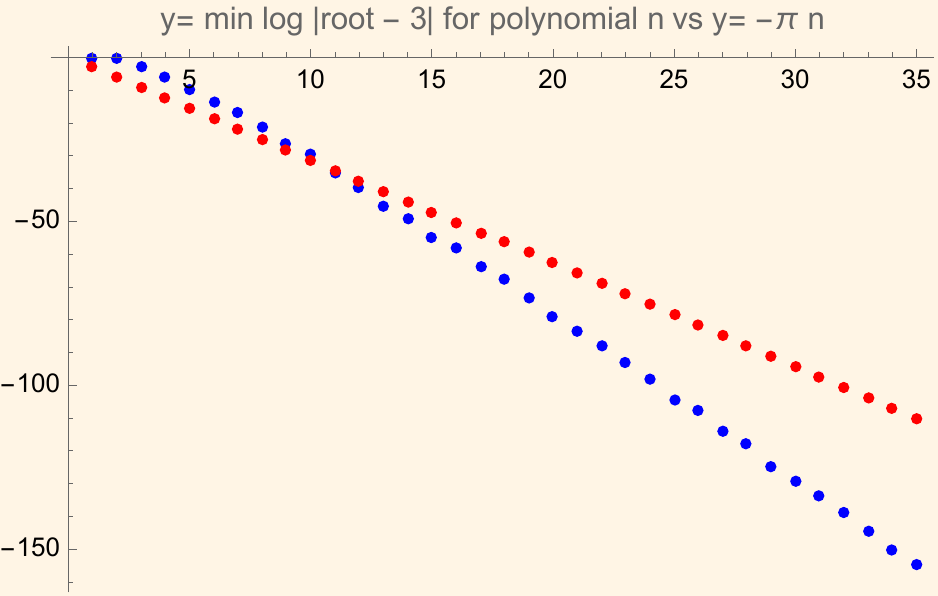}
\sc{$y$ = log(minimum distance of roots of $\overline{T}_n$ 
from $3$) in blue vs $y = -\pi n$ in red; conjectures
6 and 11.
\footnote{Notebook ``conjecture 6Laptop.nb'',
\cite{test3}.}}
\subsection{Figure 11.}
\includegraphics[scale=1]{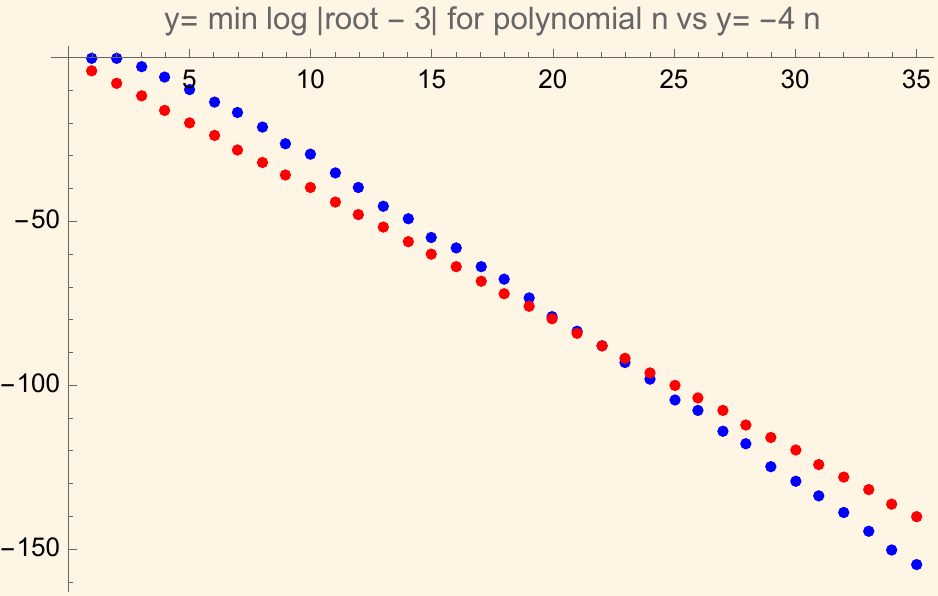}
\sc{$y$ = log(minimum distance of roots of $\overline{T}_n$ 
from $3$) in blue vs $y = -4 n$ in red; conjectures
6 and 11.
\footnote{\it ibid.\rm}}
\subsection{Figure 12.}
\includegraphics[scale=1]{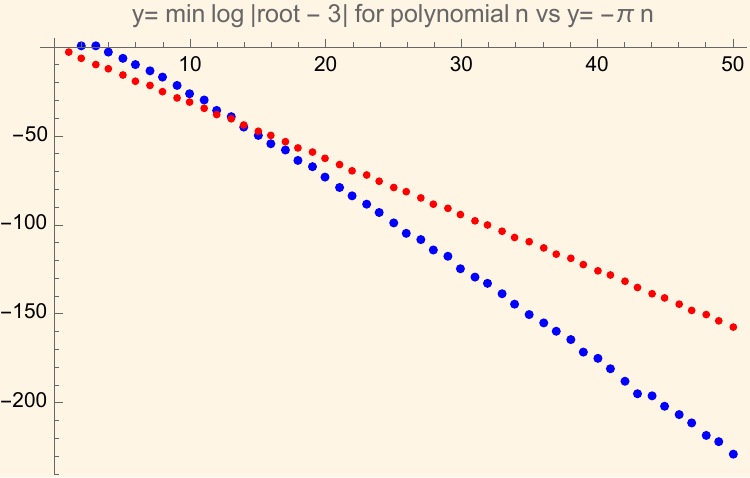}
\sc{$y$ = log(minimum distance of roots of $T_{\infty,n}$ 
from $3$) in blue vs $y = -\pi n$ in red; conjectures
7 and 11.
\footnote{Notebook ``conjecture 7.nb'',
\cite{test3}}
\subsection{Figure 13.}
\includegraphics[scale=1]{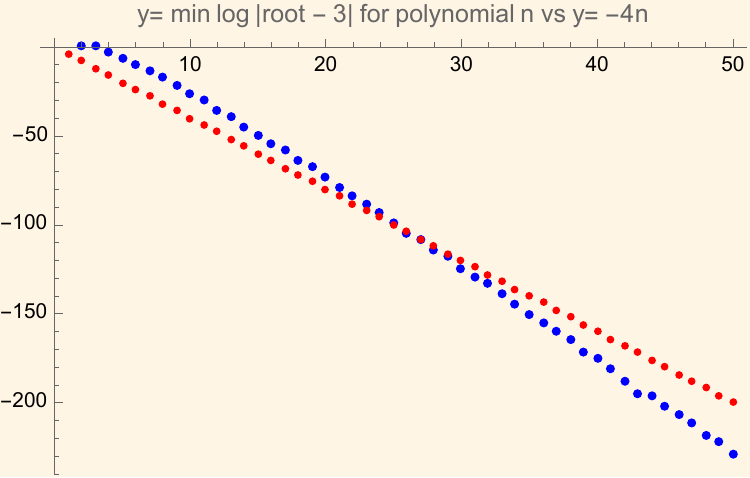}
\sc{$y$ = log(minimum distance of roots of $T_{\infty,n}$ 
from $3$) in blue vs $y = -4 n$ in red; conjectures
7 and 11.
\footnote{\it ibid.\rm}}}
\subsection{Figure 14.}
\includegraphics[scale=1]{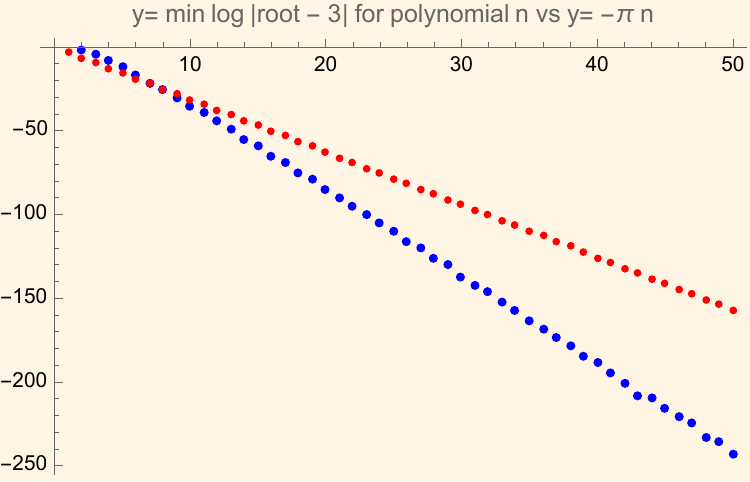}
\sc{$y$ = log(minimum distance of roots of $T^{\diamond}_n$ 
from $3$) in blue vs $y = -\pi n$ in red; conjectures
8 and 11.
\footnote{Notebook ``conjecture 8.nb'',\cite{test3}.}
\subsection{Figure 15.}
\includegraphics[scale=1]{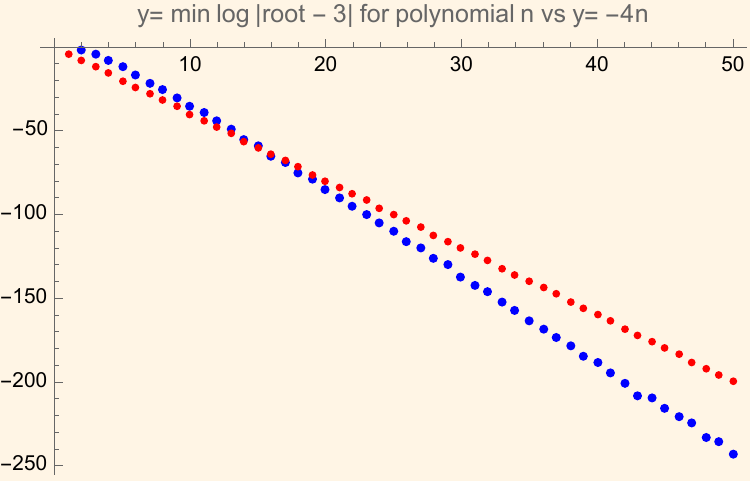}
\sc{$y$ = log(minimum distance of roots of $T^{\diamond}_n$ 
from $3$) in blue vs $y = -4 n$ in red; conjectures
8 and 11.
\footnote{\it ibid.\rm}.}
\subsection{Figure 16.}
\includegraphics[scale=1]{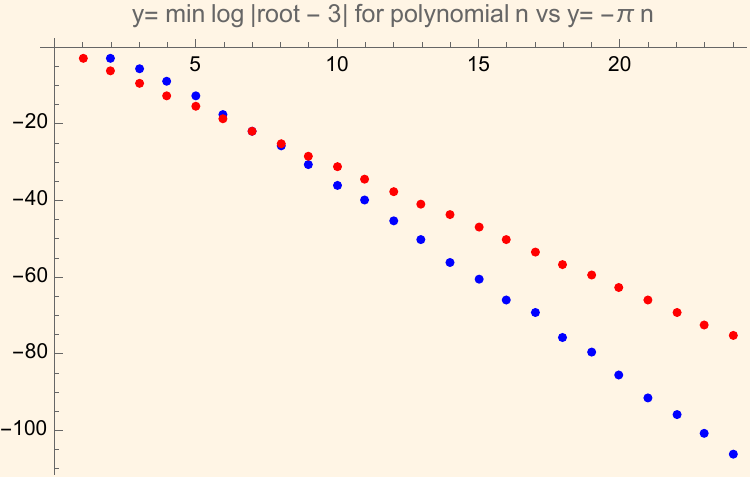}
\sc{$y$ = log(minimum distance of roots of $T^{\diamond}_{12,n}$ 
from $3$) in blue vs $y = -\pi n$ in red; conjectures
9 and 11.
\footnote{Notebook ``Conjecture 9.nb'', \cite{test3}.}}
\subsection{Figure 17.}
\includegraphics[scale=1]{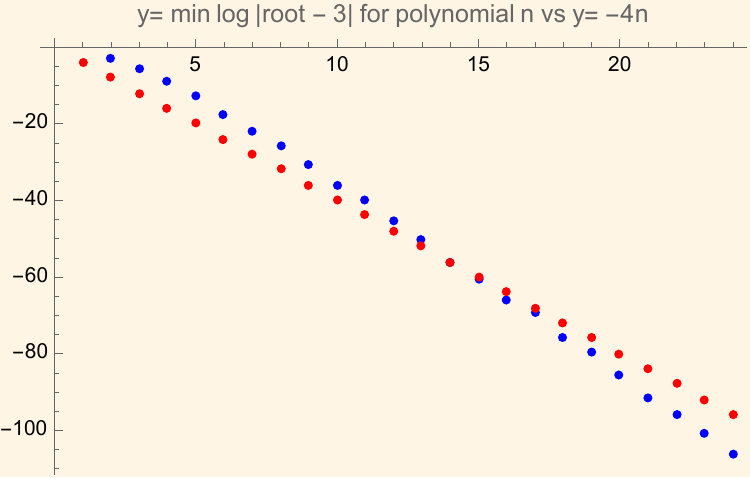}
\sc{$y$ = log(minimum distance of roots of $T^{\diamond}_{12,n}$ 
from $3$) in blue vs $y = -4 n$ in red; conjectures 9 and 11.
\footnote{\it ibid.\rm}}
\section{Tables}
\subsection[]{Table 1.}
\includegraphics[scale=.5]{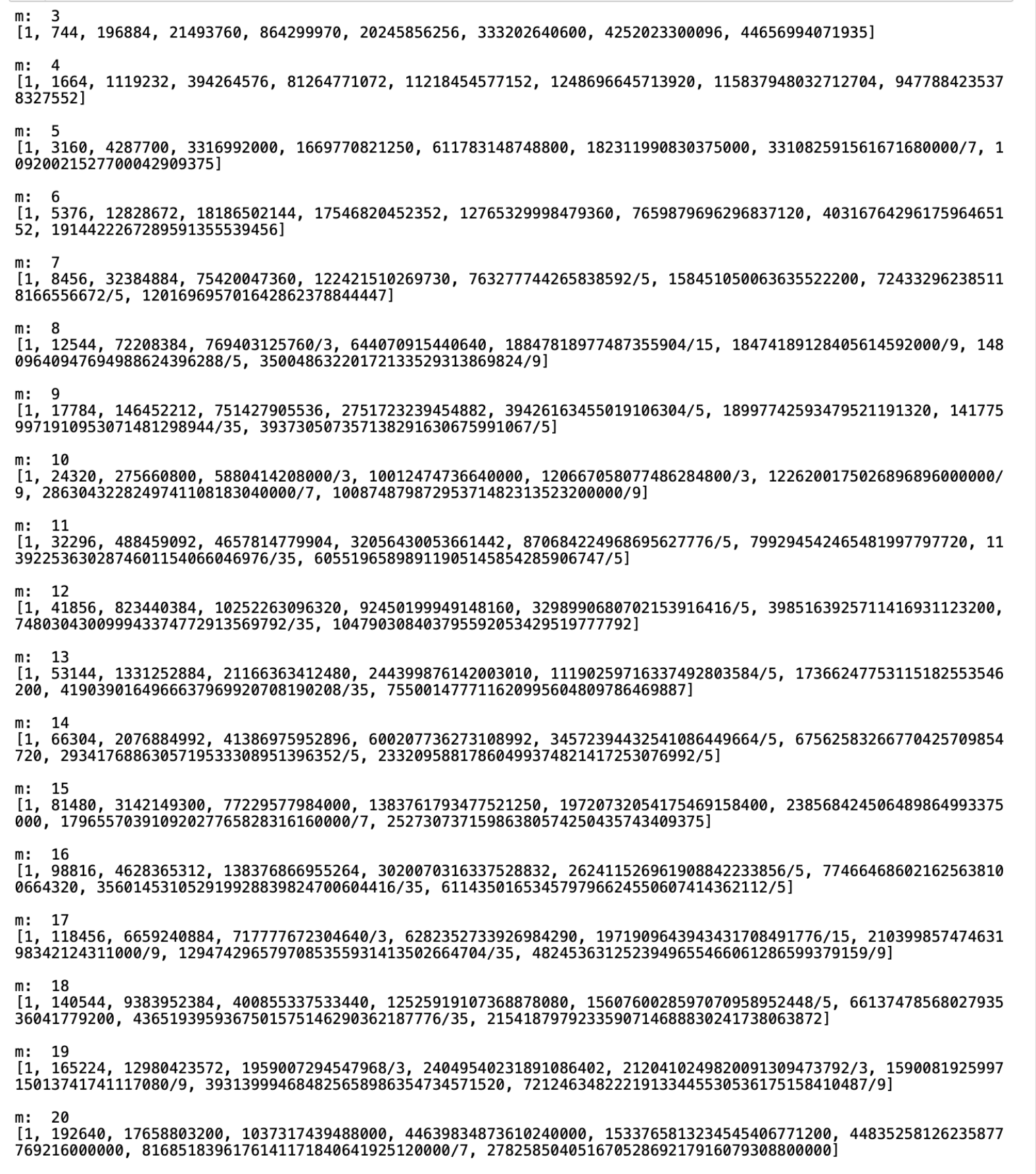}
\sc{Fourier coefficients $c_m(n)$ (conjecture 1.)
\footnote{Notebook ``conjecture 1 tables.ipynb''.}}
\newpage
\hskip -1in
\subsection[]{Table 2.}
\includegraphics[scale=.5]{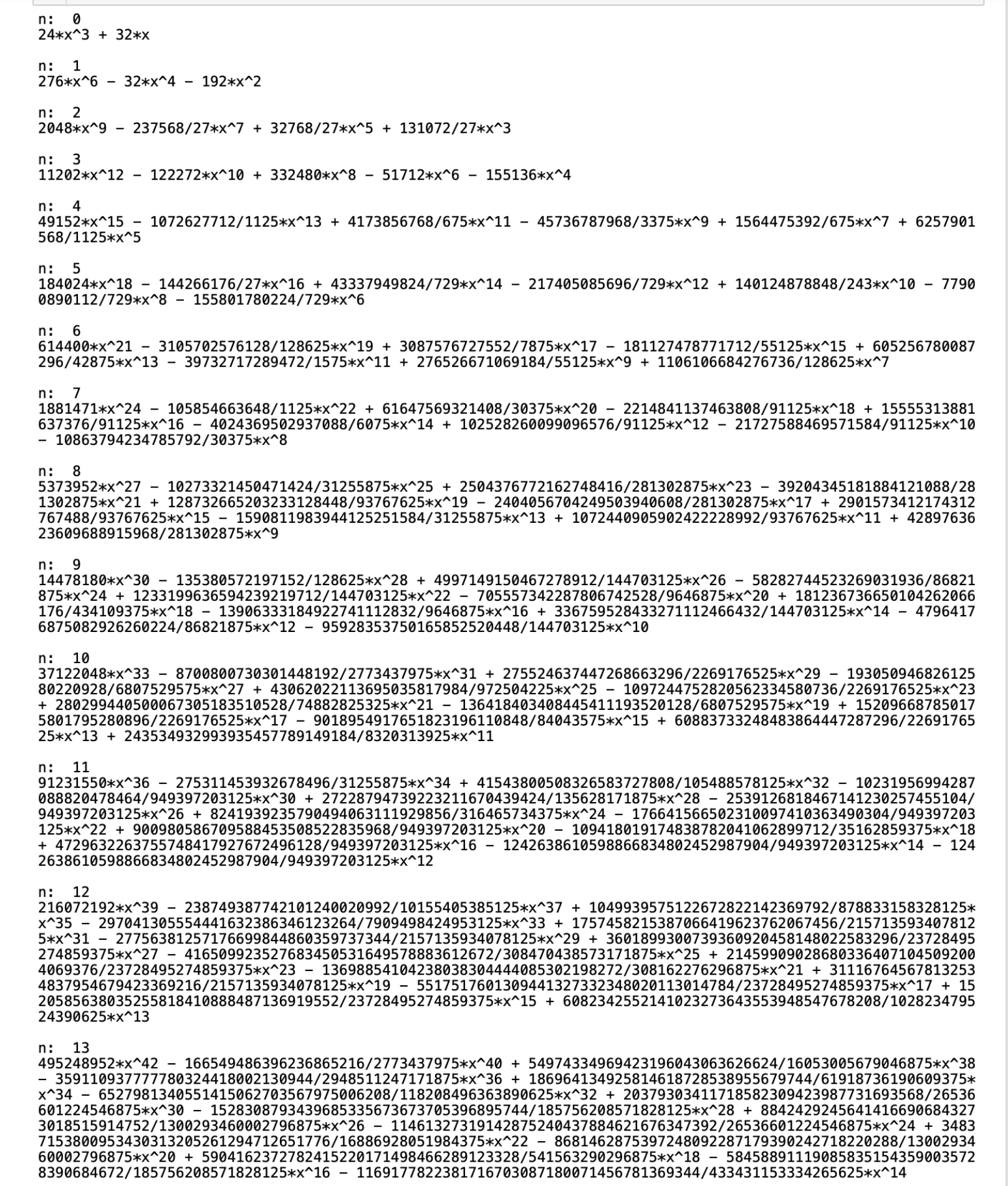}
\sc{Polynomials $C_n(x)$ (conjecture 1.)
\footnote{\it ibid.\rm}}
\newpage
\subsection[]{Table 3.}
\includegraphics[scale=.5]{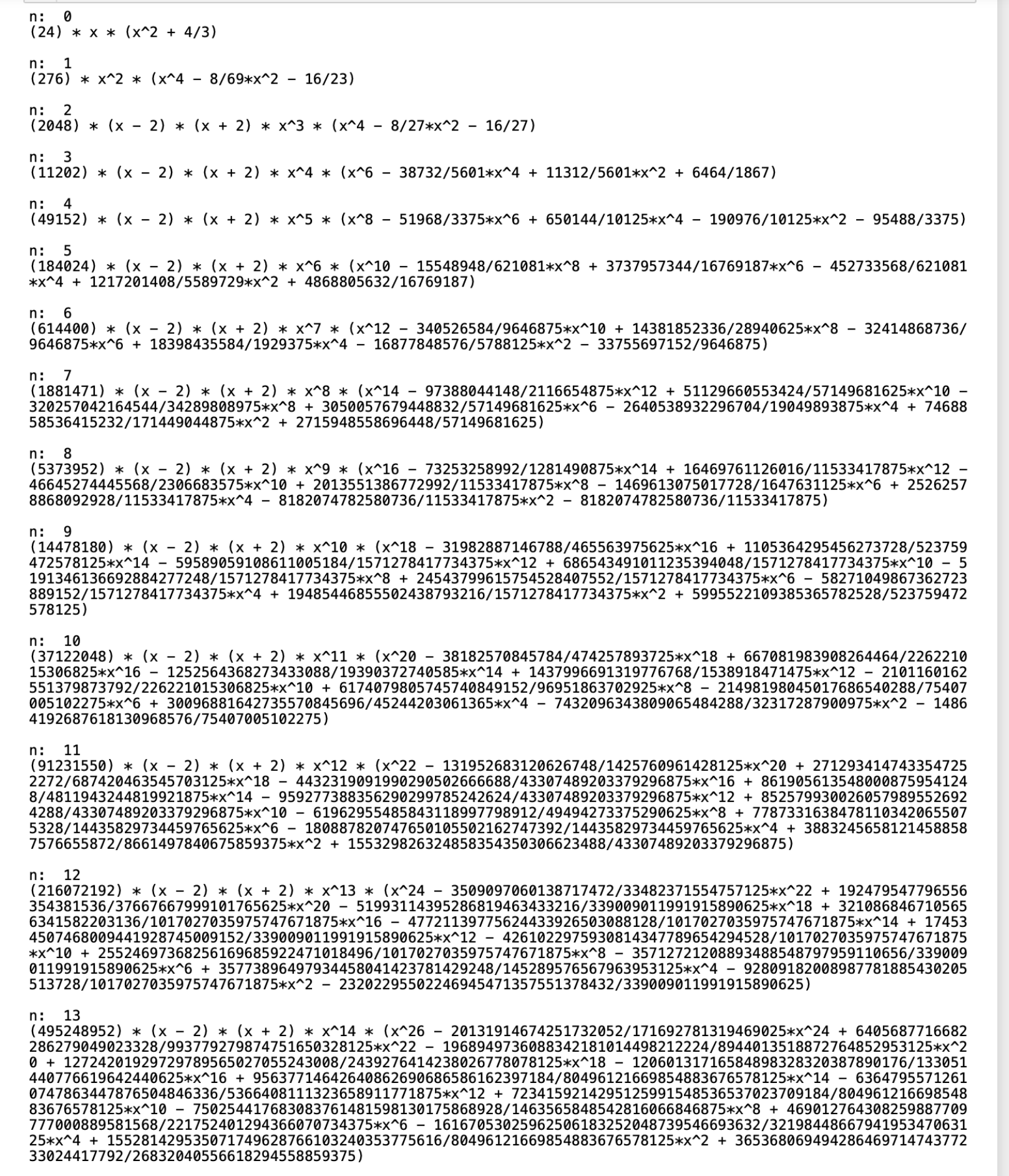}
\sc{Factored $C_n(x)$ (conjecture 1.)
\footnote{\it op. cit.\rm}}
\newpage
\subsection[]{Table 4.}
\includegraphics[scale=.3]{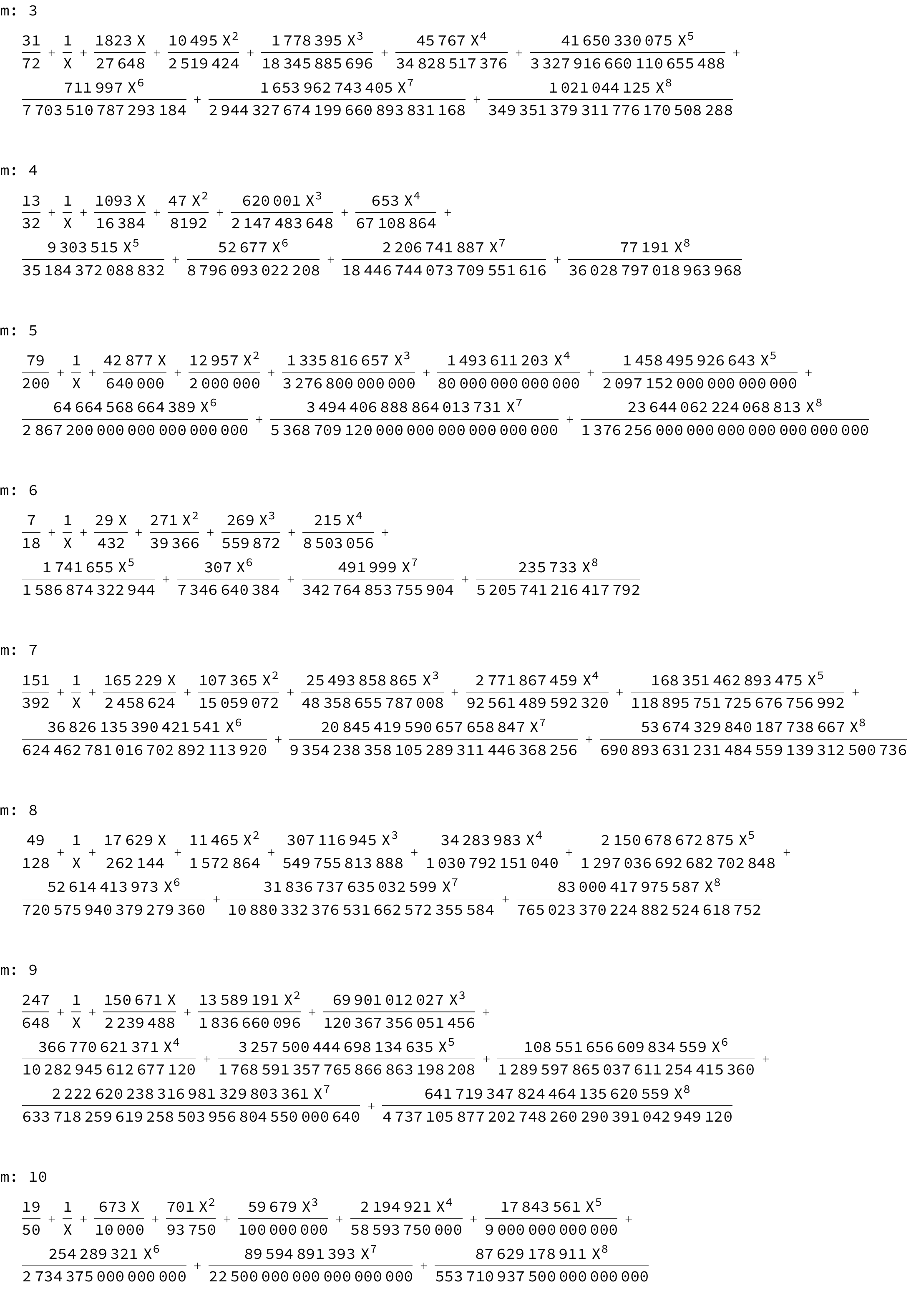}
\newline
\sc{Fourier coefficients $a_m(n)$ (conjecture 2.)
\footnote{Notebook ``conjecture 2.nb''.}}
\newpage
\hskip -1in
\subsection[]{Table 5.}
\includegraphics[scale=.3]{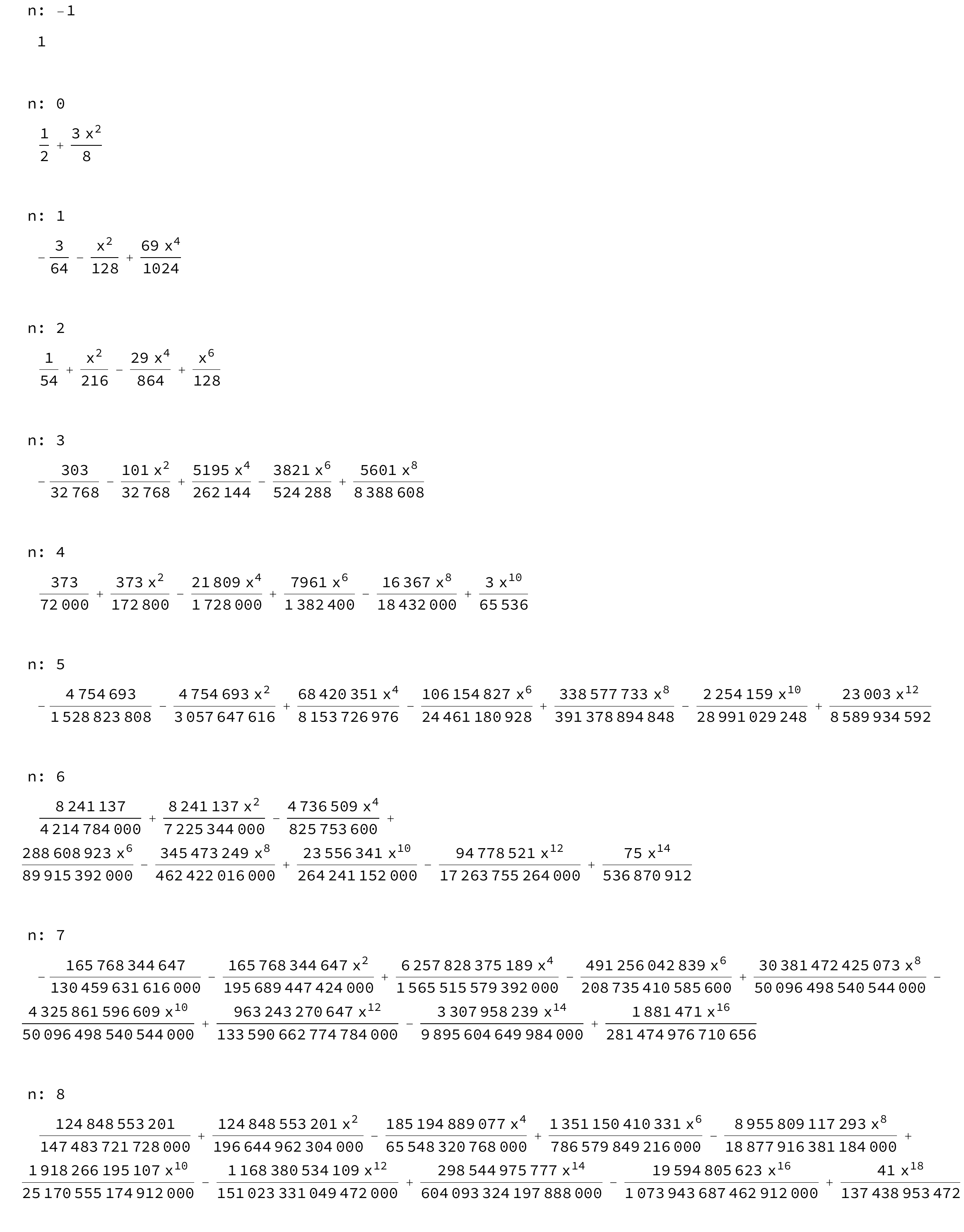}
\newline
\sc{Polynomials $A_n(x)$ (conjecture 2.)
\footnote{\it ibid.\rm}}
\newpage
\subsection[]{Table 6.}
\includegraphics[scale=.3]{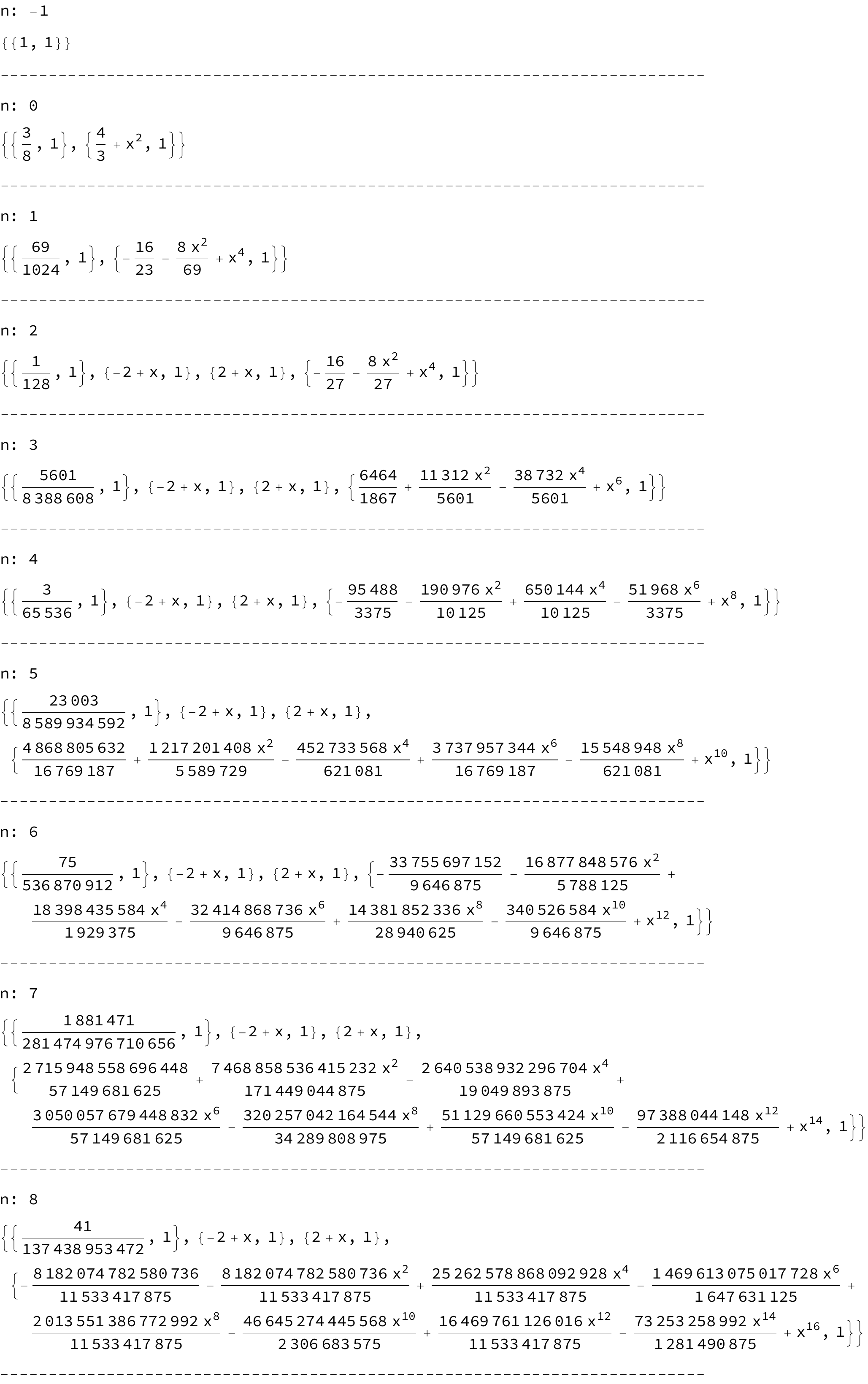}
\newline
\sc{{$A_n(x)$ factored in \it Mathematica \rm (conjecture 2.)
\footnote{\it op. cit.\rm}}
\newpage
\subsection[]{Table 7.}
\includegraphics[scale=.3]{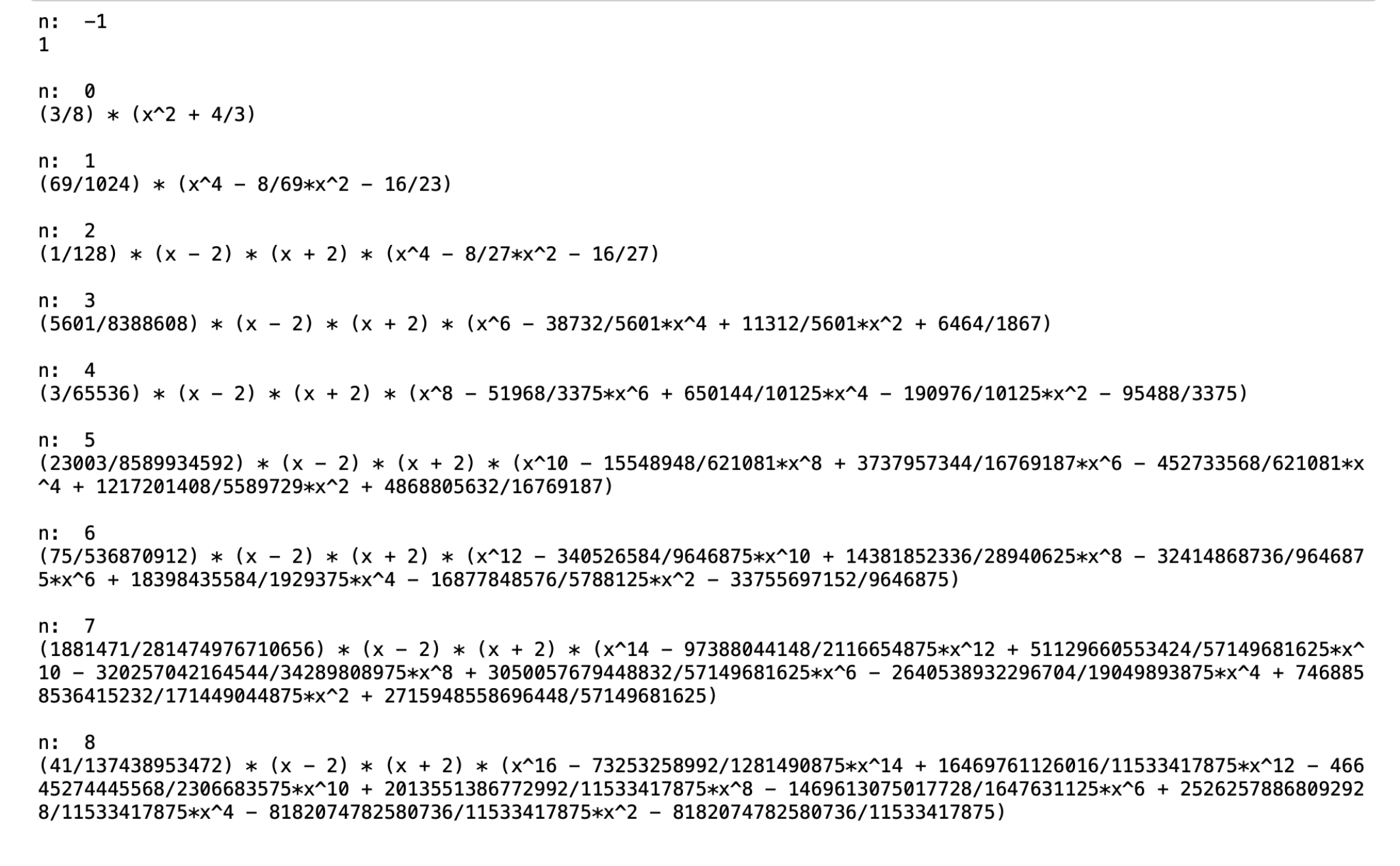}
\newline
\sc{$A_n(x)$ factored in \it SageMath \rm (conjecture 2.)
\footnote{Notebook ``conjecture 2 clause 1b.ipynb''.}}
\printbibliography
email: \tt{barrybrent``at''iphouse.com}
\end{document}